\documentclass[draft,reqno,11pt,oneside]{article}

\usepackage{amssymb}
\usepackage{english}
\usepackage{amsthm}
\usepackage{amssymb}
\usepackage{longtable}
\usepackage{amsmath}
\usepackage{ifthen}
\usepackage{upref}
\usepackage{bbold}
\usepackage{latexsym}
\usepackage[all]{xy}
\usepackage{fancyhdr}
\usepackage{makeidx}
\newfont{\chen}{cmex10 at 9pt} \newfont{\chee}{cmex10 at 10pt}
\newfont{\cheu}{cmex10 at 11pt}
\makeindex
\begin{document}
\newtheorem{lemma}{Lemma}[section]
\newtheorem*{eldemo}{Sketch of proof}
\newtheorem{satz}[lemma]{Theorem}
\newtheorem{prop}[lemma]{Proposition}
\newtheorem{defi}[lemma]{Definition}
\newtheorem{obs}[lemma]{Observation}
\newtheorem{bei}[lemma]{Example}
\newtheorem*{note}{Note}
\newtheorem{verm}[lemma]{Conjecture}
\newtheorem{kor}[lemma]{Corollary}
\renewcommand{\proofname}{Proof}
\newtheorem{bez}[lemma]{Notation}
\newtheorem{bem}[lemma]{Remark}
\newtheorem{fall}[lemma]{Case}
\newtheorem*{example}{Example}
\newtheorem{defsatz}[lemma]{Definition and Theorem}
%
\newtheorem{lemmaap}{Lemma}[section]
\newtheorem{satzap}[lemmaap]{Theorem}
\newtheorem{propap}[lemmaap]{Proposition}
\newtheorem{defiap}[lemmaap]{Definition}
\newtheorem{beiap}[lemmaap]{Example}
\newtheorem{korap}[lemmaap]{Corollary}
\newtheorem{bezap}[lemmaap]{Notation}
\newtheorem{bemap}[lemmaap]{Remark}
\def \Q{{\mathbb{Q}}}
\def \Ar{{\mathcal{A}}}
\def \G{{\mathcal{G}}}
\def \pu{{\mbox{.}}}
\def \ps{{\Phi^s}}
\def \LUC{{\mathrm{LUC}}}
\def \RUC{{\mathrm{RUC}(\G)}}
\def \RSp{{\widetilde{R}}^{(p)}}
\def \RSq{{\widetilde{R}}^{(q)}}
\def \P{{\mathcal{P}}}
\def \wwp{{\mathfrak{P}}}
\def \phii{{\overline{\Phi}}}
\def \Mun{{\mathcal{M}}}
\def \Ce{{\mathrm{C}}}
\def \Cede{{\mathbb{C}}}
\def \lpp{{L_p}}
\def \Te{{\mathcal{T}}}
\def \sis{{\mathcal{S}}}
\def \th{{\widetilde{\Theta}}}
\def \mc{{\mathcal{MC}}}
\def \ml{{\mathcal{MLUC}}}
\def \sus{{\subseteq}}
\def \mr{{\mathcal{MRUC}}}
\def \N{{\mathbb{N}}}
\def \al{{\alpha}}
\def \a{{\alpha}}
\def \suppp{{\mathrm{supp}}}
\def \cb{{\mathrm{cb}}}
\def \A{{\mathcal{A}}}
\def \llu{{\mathrm{LUC}(\G)}}
\def \B{{\mathcal{B}}}
\def \kapp{{\mathfrak{k}}}
\def \ru{{{\RUC}^*}}
\def \L{{\mathfrak{L}}}
\def \EL{{\mathcal{L}}}
\def \gip{{\Gamma_{l, p}}}
\def \UC{{{\mathrm{UC}}(\G)}}
\def \gup{{{\widetilde{\Gamma_{l, p}}}}}
\def \gurp{{\Gamma_{r, p}}}
\def \gurrp{{{\widetilde{\Gamma_{r, p}}}}}
\def \M{{\rm{M}}}
\def \en{{\mathcal{N}}}
\def \gr{{\widetilde{\Gamma_r}}}
\def \gre{{\Gamma_r}}
\def \em{{\mathcal{M}}}
\def \emli{{\mathcal{ML}_\infty}}
\def \ccc{{\mathcal{C}}}
\def \Ha{{\mathfrak{H}}}
\def \l1{{L_1(\G)}}
\def \li{{L_{\infty}(\G)}}
\def \H{{\mathcal{H}}}
\def \Qu{{\mathfrak{Q}}}
\def \cbgl{~{\stackrel{\mathrm{cb}}{=}}~}
\def \FS{{\mathcal{FS}}}
\def \gammas{{\widetilde{\gamma}}}
\def \C{{\mathcal{C}}}
\def \mr{{\mathcal{MRUC}}}
\def \Be{{\mathfrak{B}}}
\def \bar{{~|~}}
\def \Ka{{\mathcal{K}}}
\def \rr{{\mathcal{R}}}
\def \ceb{{\C\B(\B(L_2(\G)))}}
\def \alns{{\widetilde{\al_n}}} 
\def \eens{{\widetilde{e_n}}} 
\def \id{{{\mathrm{id}}}}
\def \cep{{\C\B(\B(L_p(\G)))}}
\def \dec{{\mathrm{dec}}} 
\def \ep{{\B(\B(L_p(\G)))}}
\def \Is{{\widetilde{I}}}
\def \conv{{\mathrm{conv}}}
\def \convs{{\widetilde{\mathrm{conv}}}}
\def \convsp{{\widetilde{\mathrm{conv_p}}}}
\def \es{{\B(\B(L_2(\G)))}}
\def \go{{\gip}}
\def \tt{{\overline{\otimes}}}
\def \ti{{\stackrel{\vee}{\otimes}}}
\def \oo{{{\overline{\otimes}}}}
\def \tp{{\widehat{\otimes}}}
\def \otp{{\otimes}}
\def \zet{{Z_t(\lu)}}
\def \zett{{\mathbb{Z}}} 
\makeatletter
\newcommand{\essu}
{\mathop{\operator@font ess\mbox{-}sup}}
\newcommand{\essi}
{\mathop{\operator@font ess\mbox{-}inf}}
\newcommand{\lisu}
{\mathop{\operator@font lim\,ess\mbox{-}sup}}
\newcommand{\liin}
{\mathop{\operator@font lim\,ess\mbox{-}inf}}
\makeatother
\def \ga{{\Gamma_p}}
\def \gam{{\gamma_p}}
\def \gaw{{\widetilde{\ga}}}
\def \lu{{\LUC(\G)^*}}
\def \sa{{\overline{\Gamma_p(\M(\G))}^{w*}}}
\def \Ball{{\mathrm{Ball}}}
\def \gurr{{\widetilde{\Gamma_r}}}
\def \gur{{\Gamma_r}}
\def \Ker{{\it{KERN}}}
\def \F{{\mathfrak{F}}}
\def \ad{{\mathrm{ad}}}
\def \Ceee{{\mathbb{C}}}
\def \Cee{{\mathrm{C}}}
\def \imink{{\iota_{\mathrm{minK}}}}
\def \imin{{\iota_{\mathrm{min}}}}
\def \Bild{{\it{BILD}}}
\def \Sz{{\mathcal{S}}}
\def \The{{\mathcal{T}}}
\def \cebe{{\C\B(\B(\Sz_2))}}
\def \cebes{{\C\B(\B(\Sz_2(L_2(G))))}}
\def \pr{{\Sz_2 \otimes_h \B(\Sz_2) \otimes_h \Sz_2}}
\def \su{{\rm{sup}}_{t \in \G}}
\def \mi{{\Gamma_2(\M(\G))}}
\def \Cb{{\mathcal{C}}}
\newcommand{\ens}{{\mathcal{N}}(L_p(\G))}
\newcommand{\Tee}{{\mathcal{T}}(\mathcal{H})}
\newcommand{\schutz}[1]{#1}
\def \ma{{\Gamma_2(\lu)}}
\def \Hi{{\mathcal{H}}}
\def \bs{{\B^{\sigma}(\B(\Ha))}}
\def \bsa{{\C\B(\B(\H))}}
\def \bsi{{\B^s(\B(\Ha))}}
\def \gu{{\widetilde{\Gamma_l}}}
\def \EF{{\widetilde{F}}}
\def \gi{{\Gamma_l}}
\def \als{{\widetilde{\alpha}}}
\def \betas{{\widetilde{\beta}}}
\def \gir{{\Gamma_r}}
\def \muss{{\widetilde{\mu}}}
\def \ge{{\Gamma}}
\def \gem{{\Gamma(\M(\G))}}
\def \MM{{\mathbf{M}}}
\def \un{{\ell_{\infty}^*(\G)}}
\def \R{{\mathcal{R}}}
\def \cebr{{\C\B_{\R(\G)}(\B(L_2(\G)))}}
\def \Chi{{\chi}}
\def \cebre{{\C\B_{\R(\G)}(\B(\ell_2(\G)))}}
\def \cebrs{{\C\B_{\R(\G)^{'}}^{\sigma}(\B(L_2(\G)))}}
\def \cebra{{\C\B_{\R(\G)}^{\sigma}(\B(\ell_2(\G)))}}
\def \U{{\mathfrak{U}}}
\def \K{\Ka}
\def \cebru{{\C\B_R(\B(L_2(\G)))}}
\def \lin{{\mathrm{lin}}}
\def \cebri{{\C\B_{R}^{\sigma}(\B(\H))}}
\def \Tee{{\The(\H)}}
\def \gurrn{{\gurr^0}}
\def \1{{\mathbb{1}}}
\def \gin{{\gi^0}}
\def \gun{{\gu^0}}
\newcommand{\dogl}{
\setlength{\unitlength}{0.7ex}
\linethickness {0.1ex}
~\begin{picture} (2.90 , 2)
\put (0 , 0.5) {\circle* {0.3}}
\put (0 , 1.2) {\circle* {0.3}}
\put (0.7 , 0.5) {\line (1,0) {2.1}}
\put (0.7 , 1.2) {\line (1,0) {2.1}}
\end{picture}~}
%
%
\newlength{\checklength}
\newcommand{\fns}{\footnotesize}
\newcommand{\widecheck}[1]{
\settowidth{\checklength}{$#1$}
\linethickness {0.1ex}
\setlength{\unitlength}{\checklength}
\stackrel{
\ifthenelse{ \lengthtest{ \checklength > 6em } }
{
\begin{picture}(1,0.04)
\qbezier(0,0.04)(0.3,0.03)(0.5,0.003)
\qbezier(0.5,0.003)(0.7,0.03)(1,0.04)
\qbezier(0,0.04)(0.3,0.03)(0.5,-0)
\qbezier(0.5,-0)(0.7,0.03)(1,0.04)
\qbezier(0,0.04)(0.3,0.03)(0.5,-0.005)
\qbezier(0.5,-0.005)(0.7,0.03)(1,0.04)
\end{picture}
}
{
\ifthenelse{ \lengthtest{ \checklength > 4em } }
{
\begin{picture}(1,0.05)
\qbezier(0,0.05)(0.3,0.03)(0.5,-0.02)
\qbezier(0.5,-0.02)(0.7,0.03)(1,0.05)
\qbezier(0,0.05)(0.3,0.03)(0.5,-0.015)
\qbezier(0.5,-0.015)(0.7,0.03)(1,0.05)
\qbezier(0,0.05)(0.3,0.03)(0.5,-0.01)
\qbezier(0.5,-0.01)(0.7,0.03)(1,0.05)
\end{picture}
}
{
\ifthenelse{\lengthtest{\checklength > 3em}}
{
\begin{picture}(1,0.075)
\qbezier(0,0.075)(0.3,0.05)(0.5,-0.03)
\qbezier(0.5,-0.03)(0.7,0.05)(1,0.075)
\qbezier(0,0.075)(0.3,0.05)(0.5,-0.02)
\qbezier(0.5,-0.02)(0.7,0.05)(1,0.075)
\end{picture}
}
{
\ifthenelse{\lengthtest{\checklength > 2em}}
{
\begin{picture}(1,0.1)
\qbezier(0,0.1)(0.3,0.07)(0.5,-0.04)
\qbezier(0.5,-0.04)(0.7,0.07)(1,0.1)
\qbezier(0,0.1)(0.3,0.07)(0.5,-0.02)
\qbezier(0.5,-0.02)(0.7,0.07)(1,0.1)
\end{picture}
}
{
\ifthenelse{\lengthtest{\checklength > 1em}}
{
\begin{picture}(1,0.12)
\qbezier(0,0.12)(0.3,0.1)(0.5,-0.06)
\qbezier(0.5,-0.06)(0.7,0.1)(1,0.12)
\qbezier(0,0.12)(0.3,0.1)(0.5,-0.025)
\qbezier(0.5,-0.025)(0.7,0.1)(1,0.12)
\end{picture}
}
{
\begin{picture}(1,0.15)
\qbezier(0,0.15)(0.3,0.12)(0.5,-0.1)
\qbezier(0.5,-0.1)(0.7,0.12)(1,0.15)
\qbezier(0,0.15)(0.3,0.12)(0.5,-0.04)
\qbezier(0.5,-0.04)(0.7,0.12)(1,0.15)
\end{picture}
}
}
}
}
}
}
{#1}
}
\newcommand{\swidecheck}[1]{
\settowidth{\checklength}{$_{#1}$}
\linethickness {0.1ex}
\setlength{\unitlength}{0.95\checklength}
\stackrel{
\ifthenelse{ \lengthtest{ \checklength > 6em } }
{
\begin{picture}(1,0.04)
\qbezier(0,0.04)(0.3,0.03)(0.5,0.003)
\qbezier(0.5,0.003)(0.7,0.03)(1,0.04)
\qbezier(0,0.04)(0.3,0.03)(0.5,-0)
\qbezier(0.5,-0)(0.7,0.03)(1,0.04)
\qbezier(0,0.04)(0.3,0.03)(0.5,-0.005)
\qbezier(0.5,-0.005)(0.7,0.03)(1,0.04)
\end{picture}
}
{
\ifthenelse{ \lengthtest{ \checklength > 4em } }
{
\begin{picture}(1,0.05)
\qbezier(0,0.05)(0.3,0.03)(0.5,-0.02)
\qbezier(0.5,-0.02)(0.7,0.03)(1,0.05)
\qbezier(0,0.05)(0.3,0.03)(0.5,-0.015)
\qbezier(0.5,-0.015)(0.7,0.03)(1,0.05)
\qbezier(0,0.05)(0.3,0.03)(0.5,-0.01)
\qbezier(0.5,-0.01)(0.7,0.03)(1,0.05)
\end{picture}
}
{
\ifthenelse{\lengthtest{\checklength > 3em}}
{
\begin{picture}(1,0.075)
\qbezier(0,0.075)(0.3,0.05)(0.5,-0.03)
\qbezier(0.5,-0.03)(0.7,0.05)(1,0.075)
\qbezier(0,0.075)(0.3,0.05)(0.5,-0.02)
\qbezier(0.5,-0.02)(0.7,0.05)(1,0.075)
\end{picture}
}
{
\ifthenelse{\lengthtest{\checklength > 2em}}
{
\begin{picture}(1,0.1)
\qbezier(0,0.1)(0.3,0.07)(0.5,-0.04)
\qbezier(0.5,-0.04)(0.7,0.07)(1,0.1)
\qbezier(0,0.1)(0.3,0.07)(0.5,-0.02)
\qbezier(0.5,-0.02)(0.7,0.07)(1,0.1)
\end{picture}
}
{
\ifthenelse{\lengthtest{\checklength > 1em}}
{
\begin{picture}(1,0.12)
\qbezier(0,0.12)(0.3,0.12)(0.5,-0.06)
\qbezier(0.5,-0.06)(0.7,0.12)(1,0.12)
\qbezier(0,0.12)(0.3,0.12)(0.5,-0.025)
\qbezier(0.5,-0.025)(0.7,0.12)(1,0.12)
\end{picture}
}
{
\begin{picture}(1,0.15)
\qbezier(0,0.15)(0.3,0.12)(0.5,-0.1)
\qbezier(0.5,-0.1)(0.7,0.12)(1,0.15)
\qbezier(0,0.15)(0.3,0.12)(0.5,-0.04)
\qbezier(0.5,-0.04)(0.7,0.12)(1,0.15)
\end{picture}
}
}
}
}
}
}
{#1}
}
\def \des{{\sqrt{x_j} \otimes M_{T_{\widecheck{b}}} \otimes \sqrt{x_j}}}
\def \res{{(e_j)^{\frac{1}{p}} \otimes T_{\widecheck{b}} \otimes (e_j)^{\frac{1}{q}}}}
\def \rem{{(e_j)^{\frac{1}{2}} \otimes M_{\widecheck{b}} \otimes (e_j)^{\frac{1}{2}}}}
%
%
\title{Amplification of Completely Bounded Operators \\ and
Tomiyama's Slice Maps}
\author{Matthias Neufang}
\date{}
\maketitle
\begin{abstract}
Let $(\em,\en)$ be a pair of von Neumann algebras, or 
of dual operator spaces with 
at least 
one of them having property $S_\sigma$, and let $\Phi$ be an arbitrary completely bounded mapping on $\em$. 
We\footnote{
\par
The author is currently a PIMS Postdoctoral Fellow at the
University of Alberta, Edmonton, where this work was accomplished. The support of PIMS is gratefully acknowledged.}
present an explicit construction of an amplification of
$\Phi$ 
to a completely bounded mapping on $\em \tt \en$. 
Our approach is based on the concept of slice maps as introduced by Tomiyama, and
makes use of the description of the predual of $\em \tt \en$ given by Effros and Ruan
in terms of the operator space projective tensor product (cf.\ \cite{efru}, \cite{ru}).
\par
We further discuss 
several properties of 
an amplification in connection with the investigations made in \cite{mnw}, where the special case
$\em=\B(\H)$ and $\en=\B(\Ka)$ has been considered (for Hilbert spaces $\H$ and $\K$). 
We will then 
mainly focus on various applications, such as 
a remarkable purely \textit{algebraic} characterization of $w^*$-continuity using amplifications, as well as 
a generalization of the so-called Ge--Kadison Lemma (in connection with the uniqueness problem of amplifications). 
Finally, our study 
will enable us
to show that the essential assertion of the main result in \cite{mnw} concerning completely bounded bimodule 
homomorphisms 
actually relies on
a basic property of Tomiyama's slice maps.
%
\end{abstract}

\section[Introduction]{Introduction}
\noindent
The initial aim of 
this article
is to 
present an \textit{explicit} construction of
the amplification
of an arbitrary completely bounded mapping on a von Neumann algebra $\em$ to a completely bounded mapping
on $\em \tt \en$,
where $\en$ denotes another von Neumann algebra, and to give various applications. 
\par 
Of course, an amplification of 
a completely bounded mapping on a von Neumann algebra 
is easily obtained if the latter is assumed to be {\textit{normal}} (i.e., $w^*$-$w^*$-continuous); 
cf., e.g., \cite{haa}, Lemma 1.5 (b). Our point is that we are dealing with not necessarily normal mappings 
and nevertheless even come up with an explicit formula for an amplification. 
\par 
We further
remark
that, of course,
as already established by Tomiyama, it is possible to amplify {\textit{completely positive}} mappings on
von Neumann algebras (\cite{strat}, Prop.\ 9.4). One could then think of Wittstock's decomposition
theorem to write an arbitrary completely bounded mapping on the von Neumann algebra $\em$
as a sum of four completely positive ones and apply
Tomiyama's result -- but the use of the decomposition theorem requires
$\em$ to be injective.
Moreover, we wish to stress that this procedure would be highly non-constructive:
first, the proof of Tomiyama's result uses Banach limits in an essential way, and the decomposition theorem as well
is an abstract existence result.
\par
We even go on further
to show that the same construction can actually be carried out in case $\em$
and $\en$ are only required to be dual operator spaces, where at least one of them shares property $S_\sigma$ (the latter
is satisfied, e.g., by all semidiscrete von Neumann algebras). In this abstract situation, we do not even have
a version of the above mentioned result of Tomiyama at hand.
%
\par
Furthermore, our approach yields, in particular,
a new and elegant construction of the usual amplification
of an operator $\Phi \in \C\B(\B(\H))$ to an operator in $\C\B(\B(\H \otimes_2 \Ka)) = \C\B(\B(\H) \tt \B(\Ka))$,
i.e., of the mapping
$$\Phi \mapsto \Phi^{(\infty)}$$
introduced in \cite{efki} (p.\ 265) and studied in detail in \cite{mnw} (cf.\ also \cite{bim}, p.\ 151 and Thm.\ 4.2, where
instead of $\B(\H)$, an arbitrary dual operator space is considered as the first factor).
But of course, the main interest of our approach lies
in the fact that it yields an explicit and constructive description of the amplification in a much more general setting -- where
the second factor $\B(\Ka)$ is replaced by an arbitrary von Neumann algebra or even a
dual operator space. In this case, the very definition of the mapping $\Phi \mapsto \Phi^{(\infty)}$ does not make sense,
since it is based on the representation of the elements in $\B(\H) \tt \B(\Ka)$ as infinite matrices with entries from $\B(\H)$.
\par
The crucial idea in our construction is to use the concept of slice maps as introduced by Tomiyama,
in connection with the explicit description of the predual $(\em \tt \en)_*$,
given by Effros and Ruan --
cf.\
\cite{efru}, for the case of von Neumann algebras, and \cite{ru}, for dual operator spaces with the first factor enjoying
property $S_\sigma$. In these two cases, one has a canonical complete isometry
\begin{eqnarray}
(\em \tt \en)_* \cbgl \em_* \tp \en_*, \label{projtp}
\end{eqnarray}
where the latter denotes the projective operator space tensor product.
The dichotomic nature of our statements precisely arises from this fact; at this point, we wish to emphasize that,
as shown in \cite{ru}, Cor.\ 3.7, for a dual operator space $\em$, (\ref{projtp}) holding for all dual
operator spaces $\en$, is actually equivalent to $\em$ having property $S_\sigma$.
%
\par
Hence, the
combination of operator space theory and the classical operator algebraic tool provided by
Tomiyama's slice maps enables us to show that the latter actually encodes all the
essential
properties
of an amplification.
Summarizing the major advantages of our approach to the amplification problem, we
stress the following:
\begin{itemize}
\item[(a)]
The use of Tomiyama's slice maps gives rise to
an explicit formula for the amplification of arbitrary
completely bounded mappings with a simple structure.
\item[(b)]
In particular, the approach is constructive.
\item[(c)]
In various cases, it is easier to handle amplifications using our formula,
since the construction does not involve any ($w^*$-)limits.
As we shall see,
the investigation of the amplification mapping
is thus mainly reduced to purely algebraic considerations -- in
contrast to the somehow delicate analysis used
in \cite{bim} or \cite{mnw}.
\item[(d)]
Our framework -- especially concerning the class of
objects allowed as ``second factors'' of amplifications -- is by far more general
than what can be found in the literature.
\end{itemize}
We would finally like to point out that,
in view of the above mentioned characterization of dual operator spaces fulfilling
equation (\ref{projtp}), via property $S_\sigma$ (\cite{ru}, Cor.\ 3.7),
our approach of the amplification problem
seems ``best possible'' among the approaches satisfying (a)--(d).
\par
The paper is organized as follows. -- First, we
provide the necessary terminological background from the theory of
operator algebras and
operator spaces. We further give the exact definition of what we understand by an ``amplification".
The construction of the latter in our general situation is presented in section \ref{mainconst}. 
\par
Section \ref{disc} contains a discussion of this amplification mapping under various aspects, e.g.,
relating it to results from \cite{haa} and \cite{mnw}. We derive a stability result
concerning the passage to von Neumann subalgebras or dual operator subspaces, respectively. Furthermore,
an alternative (non-constructive) description of our amplification mapping is given,
which, in particular, immediately implies that the latter preserves complete positivity.
\par
In section \ref{tow}, we
introduce an analogous concept of ``amplification in the first variable" which leads to an interesting comparison with the
usual tensor product of operators on Hilbert spaces. The new phenomenon we encounter here is that the amplifications, one in 
the first, the other in the second variable, of two completely bounded mappings $\Phi \in \C\B(\em)$ and 
$\Psi \in \C\B(\en)$ do \textit{not} commute in general. They do, though, if $\Phi$ or $\Psi$ happens to be normal (cf.\ Theorem 
\ref{commute} below). 
But our main result of this section, which may even seem a little estonishing, states that for most von Neumann algebras 
$\em$ and $\en$, the converse holds true. This means that 
the apparently weak -- and purely algebraic! -- 
condition of commutation 
of the amplification of $\Phi$ with every amplified $\Psi$ already suffices to ensure that $\Phi$ is normal (Theorem \ref{umk}). 
The proof of the latter assertion is rather technical and involves the use of countable spiral nebulas, a concept first 
introduced by Hofmeier--Wittstock in \cite{wiho}. 
\par 
In the section thereafter, 
we turn to the uniqueness problem of amplifications. Of course, an amplification which only 
meets the obvious algebraic condition is far from being unique; but we shall show that 
with some (weak) natural assumptions, one can indeed establish uniqueness. One of our two positive results 
actually yields a considerable generalization (Proposition \ref{ichgekad}) of a theorem which has proved useful in the solution of the 
splitting problem for von Neumann algebras (cf.\ \cite{gekad}, \cite{stzs}), and which is sometimes referred to as the 
Ge--Kadison Lemma. 
\par 
Finally, section \ref{sechs} is devoted to
a completely new approach to a theorem of May, Neuhardt and Wittstock regarding a
module homomorphism property of the amplified mapping. We prove a version of their result in our
context, which reveals that the essential assertion in fact relies on a fundamental and elementary
module homomorphism property shared by Tomiyama's slice maps.
%

\section{Preliminaries and basic definitions}
%
%
\par
Since
a detailed
account of the theory of operator spaces, as developed by Blecher--Paulsen, Effros--Ruan, Pisier \textit{et al.}, 
is by now available via different sources (\cite{efruan}, \cite{pis}, \cite{witt}),
as for its basic facts, we shall restrict ourselves to fix the terminology and notation we use. 
\par 
Let $X$ and $Y$ be operator spaces. We denote by $\C\B(X, Y)$ the operator space of completely bounded maps
from $X$ to $Y$, endowed with the completely bounded norm $\| \cdot \|_{\mathrm{cb}}$. We further write
$\C\B(X)$ for $\C\B(X, X)$.
If $\em$ and $\rr$ are von Neumann algebras with $\rr \subseteq \em$, then we
denote by $\C\B_{\rr} (\em)$ the space of all completely bounded
$\rr$-bimodule homomorphisms on $\em$.
If 
two operator spaces 
$X$ and $Y$ are completely isometric, then we write $X \cbgl Y$. 
\par 
An operator space $Y$ is called a \textit{dual operator space} if there is an operator space $X$ such that $Y$ is
completely isometric to $X^*$. In this case, $X$ is called an \textit{operator predual}, or \textit{predual}, for short, of
the operator space $Y$. In general, the predual of a dual operator space is not unique (up to complete isometry);
see \cite{ru}, p.\ 180, for an easy example.
\par
For a dual operator space $Y$ with a given predual $X$, we denote by $\C\B^{\sigma} (Y)$ the
subspace of $\C\B(Y)$ consisting of normal (i.e., $\sigma(Y , X)$-$\sigma(Y , X)$-continuous) mappings.
\par
If $\H$ is a Hilbert space, then, of course, every $w^*$-closed subspace $Y$ of $\B(\H)$ is a dual operator space with a
canonical predual $\B(\H)_{*}/Y_{\perp}$, where $Y_{\perp}$ denotes the preannihilator of $Y$ in $\B(\H)_*$. We denote
this canonical predual by $Y_*$.
Conversely,
if
$Y$ is a
dual operator space with a given predual $X$, there is a $w^*$-homeomorphic complete isometry of $Y$ onto a $w^*$-closed
subspace $\widetilde{Y}$ of $\B(\H)$, for some suitable Hilbert space $\H$ (see \cite{efru}, Prop.\ 5.1). Following the
terminology of \cite{efkrru}, p.\ 3, we shall call such a map a $w^*$-\textit{embedding}. In this case, $X$ is completely isometric
to the canonical predual $\widetilde{Y}_*$, and we shall identify $Y$ with $\widetilde{Y}$
(cf.\ the discussion in
\cite{ru}, p.\ 180).
\par
We denote by $\H \otimes_2 \Ka$ the Hilbert space tensor product of two Hilbert spaces $\H$ and $\Ka$.
For operator spaces $X$ and $Y$, we write $X \ti Y$ and $X \tp Y$ for the
injective and projective operator space tensor product,
respectively.
We recall that, for operator spaces $X$ and $Y$, we have a canonical complete isometry:
$$(X \tp Y)^* \cbgl \C\B(X, Y^*).$$
If $\em$ and $\en$ are von Neumann algebras, we denote as usual by $\em \tt \en$ the von Neumann algebra tensor product.
More generally, let $V^*$ and $W^*$ be dual operator spaces with given $w^*$-embeddings $V^* \subseteq \B(\H)$ and
$W^* \subseteq \B(\K)$. Then the $w^*$-spatial tensor product of $V^*$ and $W^*$, still denoted by
$V^* \tt W^*$, is defined to be the $w^*$-closure of the algebraic tensor product $V^* \otimes W^*$ in $\B(\H \otimes_2 \Ka)$.
The corresponding $w^*$-embedding of $V^* \tt W^*$ determines a predual
$(V^* \tt W^*)_*$.
Since we have
$$(V^* \tt W^*)_* \cbgl V \otimes_{\mathrm{nuc}} W$$
completely isometrically (where $\otimes_{\mathrm{nuc}}$ denotes the nuclear operator space tensor product),
the $w^*$-spatial tensor product of dual operator spaces does not depend on the given $w^*$-embedding (cf.\ \cite{efkrru}, p.\ 127).
\par
In
the following,
whenever we are speaking of ``a dual operator space", we always mean ``a dual operator space $\em$ with a given
(operator) predual, denoted by $\em_*$".
\par
We will now present, in the general framework of dual operator spaces, the construction of left and right slice maps, following
the original concept first introduced by Tomiyama for von Neumann algebras (\cite{tom}, p.\ 4;
cf.\ also \cite{kraus}, p.\ 119).
-- Let $\em$ and $\en$ be dual operator spaces. For any $\tau \in \en_*$ and any $u \in \em \tt \en$, the mapping
$$\em_* \ni \rho \mapsto \langle u, \rho \otimes \tau  \rangle$$
is a continuous linear functional on $\em_*$, hence defines an element $L_\tau (u)$ of $\em$.
As is
easily verified,
$L_\tau$ is the unique $w^*$-continuous linear map from $\em \tt \en$ to $\em$ such that
$$L_{\tau}(S \otimes T) = \langle \tau, T \rangle ~S \quad (S \in \em, ~T \in \en);$$
it will be called the left slice map associated with $\tau$. In a completely analogous fashion, we see that for every
$\rho \in \em_*$, there is a unique $w^*$-continuous linear map $R_\rho$ from $\em \tt \en$ to $\en$ such that
$$R_{\rho}(S \otimes T) = \langle \rho, S \rangle ~T \quad (S \in \em, ~T \in \en),$$
which we call the right slice map associated with $\rho$.
\par
We note that, as an immediate consequence, we have for every $\rho \in \em_*$, $\tau \in \en_*$ and
$u \in \em \tt \en$:
$$\langle L_\tau (u), \rho \rangle = \langle u, \rho \otimes \tau \rangle \quad \mathrm{and}
\quad \langle R_\rho (u), \tau \rangle = \langle u, \rho \otimes \tau \rangle.$$
Let us now briefly recall the definition of property $S_\sigma$ for dual operator spaces, as introduced by Kraus
(\cite{krtrans}, Def.\ 1.4). To this end, consider two dual operator spaces
$\em$ and $\en$ with $w^*$-embeddings $\em \subseteq \B(\H)$ and $\en \subseteq \B(\K)$.
Then the Fubini product of $\em$ and $\en$ with respect to $\B(\H)$ and $\B(\K)$ is defined through
$$\mathcal{F} (\em , \en , \B(\H) , \B(\K)) \dogl \{  u \in \B(\H \otimes_2 \K) \mid L_\tau (u) \in \em ,
~R_\rho (u) \in \en ~\mathrm{for ~all}~ \rho \in \B(\H)_* , ~\tau \in \B(\K)_* \} .$$
By the important result \cite{ru}, Prop.\ 3.3, we have
$$\mathcal{F} (\em , \en , \B(\H) , \B(\K)) \cbgl ( \em_* \widehat{\otimes} \en_* )^* ,$$
which in particular shows that the Fubini product actually does not depend on the particular choice of Hilbert spaces $\H$ and $\K$
(cf.\ also \cite{krtrans}, Remark 1.2).
Hence, we may denote the Fubini product simply by $\mathcal{F} (\em , \en)$. Now, a dual operator space $\em$
is said to have property $S_\sigma$ if
$$\mathcal{F} (\em , \en) \cbgl \em \tt \en$$
holds for every dual operator space $\en$.
\par
The question of whether every von Neumann algebra has property $S_\sigma$ has been answered in the
negative by Kraus (\cite{kraus}, Thm.\ 3.3); in fact, there exist even separably acting factors without property $S_\sigma$
(\cite{kraus}, Thm.\ 3.12).
\par
But every
semidiscrete von Neumann algebra -- and hence, every type I von Neumann algebra (\cite{efla}, Prop.\ 3.5) --
has property $S_\sigma$, as shown by Kraus in \cite{krtrans}, Thm.\ 1.9.
We finally
refer to the very interesting results
obtained by Effros, Kraus and Ruan
showing that property $S_\sigma$
is actually equivalent
to various operator space
approximation properties (cf.\ \cite{kraus}, Thm.\ 2.6; \cite{efkrru}, Thm.\ 2.4).
\par
To avoid long paraphrasing, let us
now
introduce some terminology
suitable for the statement of our results.
\begin{defi} \label{pair}
Let $\em$ and $\en$ be either
\begin{itemize}
\item[(i)] arbitrary von Neumann algebras,
\item[or]
\item[(ii)] dual operator spaces such that at least one of them has
property $S_\sigma$.
\end{itemize}
Then, in both cases, we call ($\em$, $\en$) an
${\mathrm{admissible ~pair}}$.
\end{defi}
%
%
%
If ($\em$, $\en$) is an admissible pair, then 
we have 
a canonical 
complete isometry: 
$$(\em \tt \en)_* \cbgl \em_* \tp \en_*.$$ 
%
In case $\em$ and $\en$ are both von Neumann algebras 
or dual operator spaces with $\em$ having property $S_\sigma$, 
this follows from \cite{ru}, Cor.\ 3.6 and 
Cor.\ 3.7, respectively. 
\par 
If $\em$ and $\en$ are dual operator spaces with $\en$ enjoying 
property $S_\sigma$, the above identification is seen as follows. (We 
remark that the argument below, together with \cite{ru}, Cor.\ 3.7, 
shows that 
property $S_\sigma$ of a dual operator space 
$\em$ -- resp.\ $\en$ -- is actually characterized by the equality 
$(\em \tt \en)_* \cbgl \em_* \tp \en_*$ 
holding for all dual operator spaces $\en$ -- resp.\ $\em$.) 
\par 
Let $V$ be an arbitrary operator space. 
Then, by \cite{kraus}, Thm.\ 2.6, $V^*$ has property $S_{\sigma}$ if and only if it has the 
weak${}^*$ operator approximation property (in the terminology of \cite{efkrru}). Owing to \cite{efkrru}, Thm.\ 2.4 (6), 
this is equivalent to $$V \widehat{\otimes} W \cbgl V \otimes_{\mathrm{nuc}} W$$ holding for all operator spaces $W$. 
But since both the projective and the injective operator space tensor product are symmetric, 
this in turn is equivalent to the canonical mapping $$\Phi : W \widehat{\otimes} V \longrightarrow W \stackrel{\vee}{\otimes} V$$ 
being one-to-one for all operator spaces $W$. But this is the case if and only if we have 
$$W \otimes_{\mathrm{nuc}} V \cbgl W \widehat{\otimes} V$$ for all operator spaces $W$, which finally is of course equivalent to 
$$( W^* \tt V^* )_{*} \cbgl W \widehat{\otimes} V$$ holding for all operator spaces $W$. 
\par
We
shall now make precise what we
mean by an amplification, in particular, which
topological requirements it should meet beyond the obvious algebraic ones.
\begin{defi} \label{amdef}
Let ($\em$, $\en$) be
an admissible pair. A
completely bounded linear mapping $$\Chi: \C\B(\em) \longrightarrow \C\B(\em \tt \en),$$
satisfying the
algebraic amplification condition ${\mathrm{(AAC)}}$
$$\Chi(\Phi)(S \otp T) = \Phi(S) \otp T$$
for all $\Phi \in \C\B(\em)$, $S \in \em$, $T \in \en$,
will be called
an ${\mathrm{amplification}}$, if in addition it enjoys the following properties:
\begin{itemize}
\item[(i)]
$\Chi$ is a complete isometry,
\item[(ii)]
$\Chi$ is multiplicative (hence, by (i), an isometric algebra isomorphism onto the image),
\item[(iii)]
$\Chi$ is $w^*$-$w^*$-continuous,
\item[(iv)]
$\Chi(\C\B^{\sigma}(\em)) \subseteq \C\B^{\sigma}(\em \tt \en)$, i.e., normality of the mapping is preserved.
\end{itemize}
\end{defi}
\begin{bem}
In particular, $\Chi$ is unital. For $\Chi({\mathrm{id}}_{\em})$ coincides with
${\mathrm{id}}_{\em \tt \en}$ on all elementary tensors $S \otp T \in \em \tt \en$
due to ${\mathrm{(AAC)}}$; furthermore, by (iv), it is normal,
whence we obtain
$\Chi({\mathrm{id}}_{\em}) = {\mathrm{id}}_{\em \tt \en}$ on the whole space $\em \tt \en$.
\end{bem}
In \cite{mnw} -- cf.\ also \cite{efki}, p.\ 265 --,
for every $\Phi \in \C\B(\B(\H))$, the authors obtain a mapping
$\Phi^{(\infty)} \in \C\B(\B(\H \otimes_2 \Ka)) = \C\B(\B(\H) \tt \B(\Ka))$
which obviously fulfills the algebraic amplification condition (AAC). The definition of $\Phi^{(\infty)}$ is carried out in
detail in \cite{mnw}, p.\ 284, \S 2 (there, the Hilbert space $\Ka$ is written in the form $\ell_2 (I)$ for some
suitable set $I$). We restrict ourselves to recall that, given a representation
of $u \in \B(\H \otimes_2 \Ka)$ by an infinite matrix $[u_{i , j}]_{i , j}$, where $u_{i , j} \in \B(\H)$, one has:
$$\Phi^{(\infty)} ([u_{i , j}]_{i , j}) = [\Phi (u_{i , j})]_{i , j} .$$
In fact, as is well-known, more than only (AAC) is satisfied. Indeed, we have
\begin{satz} \label{wmnek}
The mapping
\begin{eqnarray*}
\C\B(\B(\H)) &\longrightarrow& \C\B(\B(\H) \tt \B(\Ka)) \\
\Phi &\mapsto& \Phi^{(\infty)}
\end{eqnarray*}
fulfills the above requirements of an amplification.
%
\end{satz}
\begin{proof}
This follows from \cite{mnw},
p.\ 284, \S 2 and Prop.\ 2.2;
condition (iv) in
Definition \ref{amdef} is shown in \cite{hof}, Satz 3.1.
\end{proof}
We close this section by remarking that several properties of the mapping
$\Phi \mapsto \Phi^{(\infty)}$ are studied in a more general situation by Effros and Ruan in \cite{bim} (see p.\ 151 and
Thm.\ 4.2).

\section{The main construction} \label{mainconst}

%
%
We now come to our central result 
concerning the existence of an amplification in the sense of Definition \ref{amdef}
and, which is most important, presenting its explicit form. 
\begin{satz} \label{amsatz}
Let ($\em$, $\en$) be an admissible pair. Then there exists
an amplification
$$\Chi: \C\B(\em) \longrightarrow \C\B(\em \tt \en)$$
in the sense of Definition \ref{amdef}. The amplification is explicitly given by
$$\langle \Chi(\Phi) (u), \rho \otimes \tau \rangle =
\langle \Phi (L_\tau (u)), \rho \rangle ,$$
where $\Phi \in \C\B(\em)$, $u \in \em \tt \en$, $\rho \in \em_*$, $\tau \in \en_*$. -- Here, $L_\tau$ denotes
the left slice map associated with $\tau$.
\end{satz}
\begin{proof}
%
We construct a complete quotient mapping
$$\kappa: \C\B(\em \tt \en)_{*} \longrightarrow \C\B(\em)_{*}$$ in such a way that
$\Chi \dogl \kappa^{*}$ will enjoy the desired properties.
\par
We first note that
\begin{eqnarray}
\C\B(\em)_{*} \cbgl \em_* {\widehat{\otimes}} \em \label{kan1}
\end{eqnarray}
and, analogously,
\begin{eqnarray}
\C\B(\em \tt \en)_{*} \cbgl (\em \tt \en)_* {\widehat{\otimes}} (\em \tt \en) \label{kan2}
\end{eqnarray}
with completely isometric identifications.
\\
We further write
$$W: \C\B(\em \oo \en, \C\B(\en_{*}, \em)) \longrightarrow
\C\B(\en_* \tp (\em \oo \en), \em)$$
for the canonical complete isometry.
For $\tau \in \en_{*}$, let us consider
the left slice map
$$L_{\tau}: \em \oo \en \longrightarrow \em,$$ where
$$L_{\tau}(S \otimes T) = \langle \tau, T \rangle ~S$$ for $S \in \em$, $T \in \en$.
In \cite{ru}, Thm.\ 3.4 (and Cor.\ 3.6, 3.7),
Ruan constructs a complete isometric isomorphism (cf.\ also \cite{efru}, Thm.\ 3.2)
$$\theta: \em_* {\widehat{\otimes}} \en_* \longrightarrow (\em \tt \en)_*,$$
where for elementary tensors we have ($S \in \em$, $T \in \en$, $\rho \in \em_*$, $\tau \in \en_*$):
\begin{eqnarray}
\langle \theta(\rho \otimes \tau), S \otimes T \rangle = \langle \rho, S \rangle \langle \tau, T \rangle. \label{ergl}
\end{eqnarray}
Thus for all $\rho \in \em_*$ and $\tau \in \en_*$, we see that $\theta(\rho \otimes \tau)=\rho \otimes \tau$
on the algebraic tensor product $\em \otimes \en$; but since $\theta(\rho \otimes \tau), \rho \otimes \tau
\in (\em \tt \en)_*$
are normal functionals,
on the whole space $\em \tt \en$ we obtain for all $\rho \in \em_*$, $\tau \in \en_*$:
\begin{eqnarray}
\theta(\rho \otimes \tau)=\rho \otimes \tau. \label{theid}
\end{eqnarray}
Hence, in the sequel, we shall not explicitly note $\theta$ when applied on elementary tensors, except for special emphasis.
\par
We now obtain a complete isometry:
$$\theta^*: \em \tt \en \longrightarrow \C\B(\en_*, \em).$$
For each $\tau \in \en_*$, we define a linear continuous
mapping $$\varphi_{\tau}: \em {\overline{\otimes}} \en \longrightarrow \em$$
through
$$\varphi_{\tau}(u) \dogl \theta^*(u)(\tau)$$
for $u \in \em \tt \en$.
One immediately deduces from the definition that $\varphi_{\tau}$
is normal. Furthermore,
$\varphi_{\tau}$ satisfies
$$\varphi_{\tau} (S \otimes T) = \langle \tau, T \rangle ~S$$
for all $S \in \em$, $T \in \en$. -- For if $\rho \in \em_*$, we have:
\begin{eqnarray*}
\langle \varphi_{\tau}(S \otimes T), \rho \rangle &=&
\langle \theta^*(S \otimes T)(\tau), \rho \rangle \\
&=& \langle \theta^*(S \otimes T), \rho \otimes \tau \rangle \\
&=& \langle \theta(\rho \otimes \tau), S \otimes T \rangle \\
&\stackrel{(\ref{ergl})}{=}& \langle \rho, S \rangle \langle \tau, T \rangle.
\end{eqnarray*}
Hence, $\varphi_{\tau}$ coincides with the left slice map $L_\tau$ on all elementary tensors in $\em \tt \en$, and since
both mappings are normal,
we finally see that
$\varphi_{\tau}=L_{\tau}$. Thus we obtain $$\theta^*(u)(\tau) = L_{\tau} (u)$$ for all $u \in \em \tt \en$, $\tau \in \en_*$.
We note that, of course,
$$\theta^* \in \C\B(\em \tt \en, \C\B(\en_*, \em)).$$
Hence, by definition of $W$, we have:
$$W(\theta^*) \in \C\B(\en_* \tp (\em \oo \en), \em).$$
Thanks to the functorial property of the projective operator space tensor product, we now get
a completely bounded mapping
$$\id_{\em_*} \tp W(\theta^*):
\em_* \tp \en_* \tp (\em \oo \en) \longrightarrow
\em_* {\widehat{\otimes}} \em.$$
What remains to be done, is the passage from $(\em \tt \en)_* \tp (\em \tt \en)$ to
$\em_* \tp \en_* \tp (\em \oo \en)$. This is accomplished by the complete quotient mapping
$$\theta^{-1} \tp {\mathrm{id}}_{\em \tt \en}: (\em \tt \en)_* \tp (\em \tt \en)
\longrightarrow \em_* \tp \en_* \tp (\em \oo \en).$$
Now, defining
$$\kappa \dogl (\id_{\em_*} \tp W(\theta^*)) \circ (\theta^{-1} \tp {\mathrm{id}}_{\em \tt \en}): (\em \tt \en)_* \tp (\em \tt \en)
\longrightarrow \em_* {\widehat{\otimes}} \em,$$
we obtain a completely bounded mapping which has the desired properties, as we shall now prove.
\par
In view of (\ref{kan1}) and (\ref{kan2}), we see that
we have indeed constructed
a completely bounded mapping
$$\kappa: \C\B(\em \tt \en)_{*} \longrightarrow \C\B(\em)_{*},$$
and it remains to verify the properties of the mapping
$$\Chi \dogl \kappa^*: \C\B(\em) \longrightarrow \C\B(\em \tt \en).$$
\par
First, we remark the following:
\par
$(*) ~{}$
Since the mapping 
$\theta: \em_* \tp \en_* \longrightarrow (\em \tt \en)_*$ is a complete isometry,
it is readily seen that
$\theta(\em_* \otp \en_*) \otp (\em \tt \en)$ is norm dense in $(\em \tt \en)_* \tp (\em \tt \en)$.
Every element of the space $\C\B(\em \tt \en) \cbgl ((\em \tt \en)_* \tp (\em \tt \en))^*$ is thus completely determined by
its values on
$\theta(\em_* \otp \en_*) \otp (\em \tt \en)$.
\par
By construction, we see that for all $\Phi \in \C\B(\em)$, $u \in \em \tt \en$, $\rho \in \em_*$ and $\tau \in \en_*$:
\begin{eqnarray*}
\langle \Chi(\Phi), \theta(\rho \otimes \tau) \otimes u \rangle &=& \langle \Phi, \kappa(\theta(\rho \otimes \tau) \otimes u)
\rangle \\
&=& \langle \Phi, \rho \otimes W(\theta^*)(\tau \otimes u) \rangle \\
&=& \langle \Phi, \rho \otimes \theta^*(u)(\tau) \rangle \\
&=& \langle \Phi, \rho \otimes L_{\tau}(u) \rangle.
\end{eqnarray*}
Hence we have:
\begin{eqnarray*}
\langle \Chi(\Phi), \theta(\rho \otimes \tau) \otimes u \rangle =
\langle \Phi, \rho \otimes L_{\tau}(u) \rangle,
\end{eqnarray*}
and following $(*)$, by this $\Chi(\Phi)$ is completely determined. Taking into account (\ref{theid}), we thus have
that for all $\rho \in \em_*$, $\tau \in \en_*$, $u \in \em \tt \en$:
\begin{eqnarray}
\langle \Chi(\Phi), \rho \otimes \tau \otimes u \rangle =
\langle \Phi, \rho \otimes L_{\tau}(u) \rangle, \label{chi}
\end{eqnarray}
which determines $\Chi(\Phi)$ completely.
\par
Now,
let us verify
the condition (AAC) and
properties (i)--(iv).
\par
(AAC) Fix
$\Phi \in \C\B(\em)$, $S \in \em$, $T \in \en$,
$\rho \in \em_*$ and $\tau \in \en_*$.
We obtain:
\begin{eqnarray*}
\langle \rho \otimes \tau, \Chi (\Phi) (S \otimes T) \rangle &=&
\langle \Chi(\Phi), \rho \otimes \tau \otimes (S \otimes T) \rangle \\
&\stackrel{(\ref{chi})}{=}& \langle \Phi, \rho \otimes L_{\tau}(S \otimes T) \rangle \\
&=& \langle \tau, T \rangle \langle \Phi, \rho \otimes S \rangle \\
&=& \langle \rho, \Phi(S) \rangle \langle \tau, T \rangle \\
&=& \langle \rho \otimes \tau, \Phi(S) \otimes T \rangle.
\end{eqnarray*}
Hence, for all $\Phi \in \C\B(\em)$, $S \in \em$ and $T \in \en$, we have:
$$\Chi(\Phi)(S \otimes T) = \Phi(S) \otimes T,$$ as desired.
\par
(i) To prove that $\Chi$ is a complete isometry,
we show that $\kappa=\Chi_*$ is a complete quotient mapping. But this follows from the fact that
$W(\theta^*)$ is a complete quotient mapping. (For if this is established, then the mapping
$\id_{\em_*} \tp W(\theta^*) \in \C\B(\em_* \tp \en_* \tp (\em \tt \en), \em_* \tp \em)$
is a complete quotient map,
and of course this is also the case for
$\theta^{-1} \tp \id_{\em \tt \en} \in \C\B((\em \tt \en)_* \tp (\em \tt \en), \em_* \tp \en_* \tp (\em \tt \en))$.)
This in turn is easily seen. Fix $n \in \N$. Let
$S=[S_{i,j}]_{i, j} \in M_n(\em)$, $\|S\| < 1$. We are looking for an element $\varphi \in M_n(\en_* \tp (\em \tt \en))$, such that
$\| \varphi \| < 1$ and $W(\theta^*)^{(n)}(\varphi)=S$, where $W(\theta^*)^{(n)}$ denotes the $n$th amplification of the mapping
$W(\theta^*)$.
\par
Let us choose
$\tau \in
\en_*$ with $\|\tau\|=1$.
Then there exists $T \in \en$, $\|T\|=1$, such that $\langle \tau, T \rangle = 1$.
Now
$\varphi \dogl [\tau \otimes (S_{i, j} \otimes T)]_{i, j} \in M_n(\en_* \tp (\em \tt \en))$ satisfies all our requirements.
For we have:
\begin{eqnarray*}
\| \varphi \|_{M_n(\en_* \tp (\em \tt \en))} &=&
\| \tau \| ~\| [S_{i, j} \otimes T]_{i, j} \|_{M_n(\em \tt \en)} \\
&=& \| \tau \| ~\| [S_{i, j} \otimes T]_{i, j} \|_{M_n(\em \ti \en)} \\
&=& \| \tau \| ~\| [S_{i, j}]_{i, j} \|_{M_n(\em)} ~\| T \| \\
&<& 1.
\end{eqnarray*}
Furthermore, we see that:
\begin{eqnarray*}
W(\theta^*)^{(n)}(\varphi) &=& W(\theta^*)^{(n)}([\tau \otimes (S_{i, j} \otimes T)]_{i, j}) \\
&=& [W(\theta^*)(\tau \otimes (S_{i, j} \otimes T))]_{i, j} \\
&=& [\theta^*(S_{i, j} \otimes T)(\tau)]_{i, j} \\
&=& [L_{\tau} (S_{i, j} \otimes T)]_{i, j} \\
&=& [\langle \tau, T \rangle S_{i, j}]_{i, j} \\
&=& S,
\end{eqnarray*}
which establishes the claim.
\par
(ii) We show that $\Chi$ is an algebra homomorphism.
-- To this end, let $\Phi, \Psi \in \C\B(\em)$. We have to show that
$$\Chi(\Phi \Psi) = \Chi(\Phi) \Chi(\Psi),$$
as elements in $$\C\B(\em \tt \en) \cbgl (\em_* \tp \en_* \tp (\em \tt \en))^*.$$
Fix $T \in \em \tt \en$, $\rho \in \em_*$ and $\tau \in \en_*$. It suffices to show that
$$\langle \Chi(\Phi \Psi), \rho \otimes \tau \otimes T \rangle =
\langle \Chi(\Phi) \Chi(\Psi), \rho \otimes \tau \otimes T \rangle.$$
On the left side, we obtain:
\begin{eqnarray*}
\langle \Chi(\Phi \Psi), \rho \otimes \tau \otimes T \rangle &\stackrel{(\ref{chi})}{=}&
\langle \Phi \Psi, \rho \otimes L_{\tau}(T) \rangle \\
&=& \langle \Phi (\Psi(L_{\tau}(T))), \rho \rangle \\
&=& \langle \Phi, \rho \otimes \Psi(L_{\tau}(T)) \rangle.
\end{eqnarray*}
On the right, we have:
\begin{eqnarray*}
\langle \Chi(\Phi) \Chi(\Psi), \rho \otimes \tau \otimes T \rangle &=&
\langle \Chi(\Phi)[\Chi(\Psi)(T)], \rho \otimes \tau \rangle \\
&=& \langle \Chi(\Phi), \rho \otimes \tau \otimes \Chi(\Psi)(T) \rangle \\
&\stackrel{(\ref{chi})}{=}& \langle \Phi, \rho \otimes L_{\tau}(\Chi(\Psi)(T)) \rangle.
\end{eqnarray*}
Thus we have to show that
$$L_{\tau}(\Chi(\Psi)(T)) = \Psi(L_{\tau}(T)),$$ as elements in $\em$.
To this end, let $\rho \in \em_*$.
We then have:
\begin{eqnarray*}
\langle L_{\tau}(\Chi(\Psi)(T)), \rho \rangle &=& \langle \Chi(\Psi)(T), \rho \otimes \tau \rangle \\
&=& \langle \Chi(\Psi), \rho \otimes \tau \otimes T \rangle \\
&\stackrel{(\ref{chi})}{=}& \langle \Psi, \rho \otimes L_{\tau}(T) \rangle \\
&=& \langle \Psi(L_{\tau}(T)), \rho \rangle,
\end{eqnarray*}
which yields the desired equality.
\par
(iii) Being an adjoint mapping, $\Chi$ is clearly $w^*$-$w^*$-continuous.
\par
(iv) Fix $\Phi \in \C\B^{\sigma}(\em)$. We wish to prove that $\Chi(\Phi) \in \C\B^{\sigma}(\em \tt \en)$.
So let $(T_{\alpha}) \subseteq \Ball(\em \tt \en)$ be a net such that $T_{\alpha} \stackrel{w^*}{\longrightarrow} 0$.
We claim that $\Chi(\Phi)(T_{\alpha}) \longrightarrow 0$ ($\sigma(\em \tt \en, \em_* \tp \en_*)$).
It suffices to show that we have for arbitrary $\rho \in \em_*$ and $\tau \in \en_*$:
$$\langle \Chi(\Phi)(T_{\alpha}), \rho \otimes \tau \rangle \longrightarrow 0.$$
We obtain for all indices $\alpha$:
\begin{eqnarray*}
\langle \Chi(\Phi)(T_{\alpha}), \rho \otimes \tau \rangle &=& \langle \Chi(\Phi), \rho \otimes \tau \otimes T_{\alpha} \rangle \\
&\stackrel{(\ref{chi})}{=}& \langle \Phi, \rho \otimes L_{\tau}(T_{\alpha}) \rangle \\
&=& \langle \Phi(L_{\tau}(T_{\alpha})), \rho \rangle.
\end{eqnarray*}
Since $\Phi$ is normal, it thus remains to show that
$L_{\tau}(T_{\alpha})
\stackrel{w^*}{\longrightarrow} 0$; but this follows in turn from the normality of the slice map
$L_{\tau}$.
\end{proof}

\section{Discussion of the amplification mapping} \label{disc}
In the sequel, for a mapping $\Phi \in \C\B(\em)$, we shall consider amplifications with respect to different
von Neumann algebras or dual operator spaces $\en$.
We will show that all of these amplifications are compatible with each other in a very natural way.
For this purpose, we introduce the following
\begin{bez} \label{chin}
If ($\em$, $\en$) is an admissible pair, we will denote by $\Chi_{\en}$
the amplification with respect to $\en$ as constructed
in Theorem \ref{amsatz}.
\end{bez}
Taking the restriction of the mapping $\Chi_{\en}: \C\B(\em) \longrightarrow \C\B(\em \tt \en)$ to the subalgebra
$\C\B^{\sigma}(\em)$, by condition (iv) in Definition \ref{amdef}, we obtain a mapping
$\Chi_0: \C\B^{\sigma}(\em) \longrightarrow \C\B^{\sigma}(\em \tt \en)$, which of course satisfies the algebraic condition
(AAC) on $\C\B^{\sigma}(\em)$.
But this is the case as well for the amplification
$$H: \C\B^{\sigma}(\em) \longrightarrow \C\B^{\sigma}(\em \tt \en)$$
considered, for von Neumann algebras $\em$ and $\en$, by
de Canni\`ere and Haagerup in \cite{haa}, Lemma 1.5 (b).
Note that $H$
is completely determined by
(AAC), since the image of $H$ consists of normal mappings.
Hence we deduce that $\Chi_0=H$, thus showing that $\Chi_{\en}$ is indeed a ($w^*$-$w^*$-continuous) extension of $H$
to $\C\B(\em)$, where $\em$ and $\en$ are von Neumann algebras (cf.\ also \cite{kraus}, p.\ 123, for a
discussion of the mapping
$H$ in a more general setting). 
\begin{bem} 
Theorem \ref{amsatz} 
implies in particular (cf.\ (iv) in Definition \ref{amdef}) that 
whenever $(\em , \en)$ is an admissible pair, then for every $\Phi \in \C\B^\sigma (\em)$, 
the mapping $\Chi_\en (\Phi)$ is the unique normal map in $\C\B^\sigma (\em \tt \en)$ which satisfies (AAC). 
This entails that the last assumption made in \cite{krtrans}, Prop.\ 1.19 is always fulfilled (and even true in greater 
generality than actually needed there), which besides slightly simplifies the proof of \cite{krtrans}, Thm.\ 2.2. 
\end{bem} 
Let us now briefly return to the special case where $\em=\B(\H)$ and $\en=\B(\Ka)$.
We already
noted in Theorem \ref{wmnek} that the mapping
$\Phi \mapsto \Phi^{(\infty)}$, where $\Phi \in \C\B(\B(\H))$, is an amplification
in the sense of Definition \ref{amdef}. Since $\B(\H)$ trivially
is an injective factor, we conclude by Proposition \ref{eindam} that
\begin{eqnarray}
\Chi_{\B(\Ka)} (\Phi) = \Phi^{(\infty)} \label{glinf}
\end{eqnarray}
for all $\Phi \in \C\B(\B(\H))$.
Thus, restricted to that special case,
Theorem \ref{amsatz} provides
a new construction of the mapping $\Phi \mapsto \Phi^{(\infty)}$ and for the first time shows the intimate relation to
Tomiyama's slice
maps.
\par
We finally stress that
our construction of $\Chi_{\B(\Ka)} (\Phi)$ is coordinate free -- in contrast to the
definition of $\Phi^{(\infty)}$ in \cite{mnw}:
since
our approach does not rely on an
explicit representation of operators in $\B(\H) \tt \B(\Ka)$
as infinite matrices, we do not have to choose a basis in the Hilbert space $\Ka$, whereas in \cite{mnw}, the latter is
represented as $\ell_2(I)$ for some suitable $I$. The equality
$$\Chi_{\B(\Ka)} (\Phi) = \Phi^{(\infty)} \quad (\Phi \in \C\B(\B(\H))),$$
obtained above, now immediately implies that
the definition of the mapping $\Phi^{(\infty)}$ in \cite{mnw} does not depend on the particular choice of a basis. (This fact
of course can also be obtained intrinsically from the construction of $\Phi^{(\infty)}$, but it comes for free in our approach.)
\par
We now come to the result announced above
which establishes the compatibility of our amplification with respect to different
von Neumann algebras or dual operator spaces $\en$, respectively.
\begin{satz} \label{kompa}
Let ($\em$, $\en$) be an admissible pair. Let further $\en_0 \subseteq \en$ be either a von Neumann subalgebra
or a dual operator subspace of $\en$, respectively (according to case (i) or (ii) of Definition \ref{pair}).
Then, for all
$\Phi \in \C\B(\em)$, we have:
$$\Chi_{\en}(\Phi) |_{\em \tt \en_0} = \Chi_{\en_0}(\Phi).$$
\end{satz}
\begin{proof}
Denote by $\iota: \em \tt \en_0 \hookrightarrow \em \tt \en$ the canonical embedding, and write
$q: \en_* \twoheadrightarrow (\en_0)_*$
for the
restriction map. Obviously, the pre-adjoint mapping $\iota_*: \em_* \tp \en_* \twoheadrightarrow
\em_* \tp (\en_0)_*$ satisfies
$\iota_* = \id_{\em_*} \tp q$.
\par
We
shall prove that the following diagram commutes:
$$\xymatrix{
\em \tt \en \ar[r]_{\Chi_{\en}(\Phi)} & \em \tt \en \\
\em \tt \en_0 \ar@{^{(}->}[u]_{\iota} \ar[r]^{\Chi_{\en_0}(\Phi)} & \em \tt \en_0 \ar@{^{(}->}[u]^{\iota}
}
$$
This yields the desired equality (and shows at the same time that
$\Chi_{\en}(\Phi) |_{\em \tt \en_0}$
indeed
leaves the space $\em \tt \en_0$ invariant).
\par
To this end, fix $u \in \em \tt \en_0$. We have to show that
$$\Chi_{\en}(\Phi)(\iota(u)) = \iota(\Chi_{\en_0}(\Phi)(u)).$$
Let $\rho \in \em_*$ and $\tau \in \en_*$. It is sufficient to verify that
$$\langle \Chi_{\en}(\Phi)(\iota(u)), \rho \otp \tau \rangle =
\langle \iota(\Chi_{\en_0}(\Phi)(u)), \rho \otp \tau \rangle.$$
For the left-hand expression, we get:
$$\langle \Chi_{\en}(\Phi)(\iota(u)), \rho \otp \tau \rangle =
\langle \Chi_{\en}(\Phi), \rho \otp \tau \otimes
\iota(u) \rangle = \langle \Phi, \rho \otp L_{\tau}(\iota(u)) \rangle.$$
On the right-hand side, we find:
\begin{eqnarray*}
\langle \iota(\Chi_{\en_0}(\Phi)(u)), \rho \otp \tau \rangle &=&
\langle \Chi_{\en_0}(\Phi)(u), \rho \otp q(\tau) \rangle \\ &=&
\langle \Chi_{\en_0}(\Phi), \rho \otp q(\tau) \otp u \rangle \\ &=&
\langle \Phi, \rho \otimes L_{q(\tau)}(u) \rangle.
\end{eqnarray*}
It thus remains to be shown that
$$L_{\tau}(\iota(u)) = L_{q(\tau)}(u).$$
But for $\varphi \in \em_*$, we obtain:
\begin{eqnarray*}
\langle L_{\tau}(\iota(u)), \varphi \rangle &=&
\langle \iota(u), \varphi \otp \tau \rangle \\ &=&
\langle u, \varphi \otp q(\tau) \rangle \\ &=&
\langle L_{q(\tau)}(u), \varphi \rangle,
\end{eqnarray*}
which yields the claim.
\end{proof}
%
%
We shall now present another very natural, though highly non-constructive,
approach to the amplification problem, and discuss its
relation to
our concept.
-- Let ($\em$, $\en$) be an admissible pair, with $\em \subseteq \B(\H)$ and $\en \subseteq \B(\K)$.
An element $\Phi \in \C\B(\em)$ may
be considered as a completely bounded mapping from $\em$ with
values in $\B(\H)$. Hence we can assign to $\Phi$ a
Wittstock-Hahn-Banach extension $\widetilde{\Phi} \in \C\B(\B(\H))$.
Then the mapping
$$\widetilde{\Phi}^{(\infty)} \mid_{\em \tt \en}: \em \tt \en \longrightarrow \B(\H) \tt \B(\K)$$
is completely bounded and of course satisfies the algebraic amplification condition (AAC).
As we shall now see, this (non-constructive) abstract procedure yields indeed a mapping
which takes
values in $\em \tt \en$ and
does not depend on the choice of the Wittstock-Hahn-Banach extension --
namely, $\widetilde{\Phi}^{(\infty)} \mid_{\em \tt \en}$ is nothing but $\Chi_\en (\Phi)$. This will
be proved in a similar fashion as Theorem \ref{kompa}.
\begin{satz} \label{ichuwitt}
Let ($\em$, $\en$) be an admissible pair, where $\em \subseteq \B(\H)$ and $\en \subseteq \B(\K)$.
Let further $\Phi \in \C\B(\em)$. Then for an arbitrary
Wittstock-Hahn-Banach extension $\widetilde{\Phi} \in \C\B(\B(\H))$ obtained as above, we have:
$$\widetilde{\Phi}^{(\infty)} \mid_{\em \tt \en} = \Chi_\en (\Phi).$$
\end{satz}
\begin{proof}
We write $\iota: \em \tt \en \hookrightarrow \B(\H) \tt \B(\K)$ for the canonical embedding. It is sufficient to
show that the following diagram commutes:
$$\xymatrix{
\B(\H) \tt \B(\K) \ar[r]_{\widetilde{\Phi}^{(\infty)}} & \B(\H) \tt \B(\K) \\
\em \tt \en \ar@{^{(}->}[u]_{\iota} \ar[r]^{\Chi_{\en}(\Phi)} & \em \tt \en \ar@{^{(}->}[u]^{\iota}
}
$$
Fix $u \in \em \tt \en$, $\rho \in \B(\H)_*$, $\tau \in \B(\K)_*$.
We have only to show that
$$\langle \iota [ \Chi_\en (\Phi) (u)  ], \rho \otimes \tau \rangle =
\langle \widetilde{\Phi}^{(\infty)} (\iota (u)), \rho \otimes \tau \rangle ,$$
which by equation (\ref{glinf}) is equivalent to
\begin{eqnarray}
\langle \iota [ \Chi_\en (\Phi) (u)  ], \rho \otimes \tau \rangle =
\langle \Chi_{\B(\K)} (\widetilde{\Phi})  (\iota (u)), \rho \otimes \tau \rangle . \label{anfan}
\end{eqnarray}
On the left-hand side, we obtain:
\begin{eqnarray*}
\langle \iota [ \Chi_\en (\Phi) (u)  ], \rho \otimes \tau \rangle &=&
\langle \Chi_\en (\Phi) (u) , \rho \mid_\em \otimes~ \tau \mid_\en \rangle \\
&=& \langle \Phi, \rho \mid_\em \otimes~ L_{\tau \mid_\en} (u)  \rangle .
\end{eqnarray*}
Now, writing $\iota' : \em \hookrightarrow \B(\H)$ for the canonical embedding, we note that
\begin{eqnarray*}
\langle L_\tau (\iota (u)), \rho \rangle &=& \langle \iota (u), \rho \otimes \tau \rangle \\
&=& \langle u, \rho \mid_\em \otimes~ \tau \mid_\en \rangle \\
&=& \langle L_{\tau \mid_{\en}} (u), \rho \mid_\em \rangle \\
&=& \langle \iota' (L_{\tau \mid_{\en}} (u)), \rho \rangle ,
\end{eqnarray*}
whence we have:
$$L_\tau (\iota (u)) = \iota' (L_{\tau \mid_{\en}} (u)).$$
Thus, we see that
\begin{eqnarray}
\widetilde{\Phi} [L_\tau (\iota (u))] = \widetilde{\Phi} [\iota' (L_{\tau \mid_{\en}} (u))]
= \iota' [\Phi (L_{\tau \mid_{\en}} (u))] . \label{letzt}
\end{eqnarray}
Now we find for the right-hand side term of equation (\ref{anfan}):
\begin{eqnarray*}
\langle \Chi_{\B(\K)} (\widetilde{\Phi})  (\iota (u)), \rho \otimes \tau \rangle &=&
\langle \widetilde{\Phi} , \rho \otimes L_\tau (\iota (u))  \rangle \\
&=& \langle \widetilde{\Phi} [L_\tau (\iota (u))] , \rho \rangle \\
&\stackrel{(\ref{letzt})}{=}& \langle \iota' [\Phi (L_{\tau \mid_{\en}} (u))] , \rho \rangle \\
&=& \langle \Phi (L_{\tau \mid_{\en}} (u)) , \rho \mid_{\em}   \rangle \\
&=& \langle \Phi, \rho \mid_\em \otimes~ L_{\tau \mid_\en} (u)  \rangle ,
\end{eqnarray*}
whence we deduce the desired equality.
\end{proof}
\begin{bem}
If $\em \subseteq \B(\H)$, $\en \subseteq \B(\K)$ are von Neumann algebras and $\Phi \in \C\P(\em)$ is completely positive, then,
using an arbitrary Arveson extension of $\Phi$ to a completely positive map $\widetilde{\Phi} \in \C\P(\B(\H))$, we obtain by
the same reasoning as above that $\widetilde{\Phi}^{(\infty)} \mid_{\em \tt \en} = \Chi_\en (\Phi)$. But it is
easy to see from the definition of the mapping $\widetilde{\Phi} \mapsto \widetilde{\Phi}^{(\infty)}$
(cf.\ \cite{mnw}, p.\ 284, \S 2)
that $\widetilde{\Phi}^{(\infty)}$ still is completely positive, whence $\Chi_\en (\Phi)$ also is.
This shows that our amplification $\Chi_\en$ preserves complete positivity.
\end{bem}

\section{An algebraic characterization of 
normality} \label{tow}
In the following, let still ($\em$, $\en$) be an admissible pair.
Using Tomiyama's left slice map, for every completely bounded mapping
$\Phi \in \C\B(\em)$,
we have constructed an amplification
$\Chi_{\en}(\Phi) =: \Phi \otp \id_{\en} \in \C\B(\em \tt \en)$.
In an analogous fashion, replacing in the above construction Tomiyama's left by the appropriate right slice
maps,
we obtain -- {\textit{mutatis mutandis}} -- an amplification of completely bounded mappings
$\Psi \in \C\B(\en)$
``in the left variable".
This yields an amplification
$\id_{\em} \otp \Psi \in \C\B(\em \tt \en)$, which is completely determined by the equation:
\begin{eqnarray}
\langle \id_{\em} \otp \Psi, \rho \otp \tau \otp u \rangle = \langle \Psi, \tau \otp R_\rho(u) \rangle, \label{defr}
\end{eqnarray}
where $u \in \em \tt \en$, $\rho \in \em_*$ and $\tau \in \en_*$ are arbitrary.
Of course, $\id_{\em} \otp \Psi$ shares analogous properties to those of
$\Phi \otp \id_{\en}$.
\par
It is natural to ask for some relation between the operators in
$\C\B(\em \tt \en)$ which arise
by the two different kinds of amplification.
Let us briefly compare our situation with the classical setting of amplification of
operators on Hilbert spaces $\H$ and $\Ka$.
For operators $S \in \B(\H)$ and $T \in \B(\Ka)$, we consider the amplifications $S \otp \id_\Ka$ and
$\id_\H \otp T$ in $\B(\H) \tt \B(\Ka) = \B(\H \otp_2 \Ka)$. Then it is evident that we have
\begin{eqnarray}
(S \otp \id_\Ka) ~(\id_\H \otp T) = (\id_\H \otp T) ~(S \otp \id_\Ka), \label{commu}
\end{eqnarray}
and the tensor product of $S$ and $T$ is defined precisely to be this operator.
\par
Back in our situation,
for $\Phi \in \C\B(\em)$ and $\Psi \in \C\B(\en)$,
the corresponding amplifications $\Phi \otp \id_\en$ and $\id_\em \otp \Psi$ will 
commute provided at least one of the involved mappings $\Phi$ or $\Psi$ is normal -- 
but not 
in general. 
In fact, we will even show that for most von Neumann algebras $\em$ and $\en$, 
$\Phi \in \C\B (\em)$ satisfying this innocent looking commutation relation for all $\Psi \in \C\B(\en)$ 
forces $\Phi$ to be normal! 
This phenomenon, which may be surprising at first glance, 
is due to the fact that, in contrast to the multiplication \textit{in} a von Neumann algebra $\rr$ (in 
our case $\rr=\em \tt \en$),
the multiplication in $\C\B(\rr)$ is not $w^*$-continuous in both variables, but only in the left one. In fact,
the subset of $\C\B(\rr) \times \C\B(\rr)$ on which the product is $w^*$-continuous is
precisely $\C\B(\rr) \times \C\B^{\sigma}(\rr)$.
This is reflected in the following result which parallels relation (\ref{commu}) in our context.
\begin{satz} \label{commute}
Let ($\em$, $\en$) be an admissible pair.
If $\Psi \in \C\B^\sigma(\en)$ is a normal mapping, then we
have for every $\Phi \in \C\B(\em)$:
$$(\Phi \otp \id_\en) ~(\id_\em \otp \Psi) = (\id_\em \otp \Psi) ~(\Phi \otp \id_\en).$$
\end{satz}
%
%
\begin{proof}
Fix $u \in \em \tt \en$, $\rho \in \em_*$ and $\tau \in \en_*$.
We first
remark that
\begin{eqnarray}
\langle \id_{\em} \otp \Psi, \rho \otp \tau \otp u \rangle &\stackrel{(\ref{defr})}{=}& \langle \Psi, \tau \otp R_\rho(u) \rangle
\nonumber\\
&=& \langle \Psi(R_\rho(u)), \tau \rangle \nonumber\\
&=& \langle R_\rho(u), \Psi_*(\tau) \rangle \nonumber\\
&=& \langle u, \rho \otp \Psi_*(\tau) \rangle, \label{urops}
\end{eqnarray}
where $\Psi_*$ denotes the pre-adjoint mapping of $\Psi$.
\par
To establish the claim, it suffices to prove that:
$$\langle ( \Phi \otp \id_\en ) ~( \id_\em \otp \Psi ), \rho \otp \tau \otp u \rangle =
\langle ( \id_\em \otp \Psi ) ~( \Phi \otp \id_\en ), \rho \otp \tau \otp u \rangle.$$
The left side takes on the following form:
\begin{eqnarray*}
\langle ( \Phi \otp \id_\en ) ~( \id_\em \otp \Psi ), \rho \otp \tau \otp u \rangle &=&
\langle \Chi_\en(\Phi), \rho \otp \tau \otp ( \id_\em \otp \Psi )(u) \rangle \\
&=& \langle \Phi, \rho \otp L_\tau [( \id_\em \otp \Psi )(u)] \rangle,
\end{eqnarray*}
whereas on the right side, we obtain:
\begin{eqnarray*}
\langle ( \id_\em \otp \Psi ) ~( \Phi \otp \id_\en ), \rho \otp \tau \otp u \rangle &=&
\langle ( \id_\em \otp \Psi ) ~( \Chi_\en(\Phi) (u) ), \rho \otp \tau \rangle \\
&\stackrel{(\ref{urops})}{=}& \langle \Chi_\en(\Phi) (u), \rho \otp \Psi_*(\tau) \rangle \\
&=& \langle \Chi_\en(\Phi), \rho \otp \Psi_*(\tau) \otp u \rangle \\
&=& \langle \Phi, \rho \otp L_{\Psi_*(\tau)} (u) \rangle.
\end{eqnarray*}
Hence, it remains to be shown that
$$L_\tau [( \id_\em \otp \Psi )(u)] = L_{\Psi_*(\tau)} (u).$$
But for arbitrary $\rho' \in \em_*$, we get:
\begin{eqnarray*}
\langle L_\tau [( \id_\em \otp \Psi )(u)], \rho' \rangle &=& \langle ( \id_\em \otp \Psi )(u), \rho' \otp \tau \rangle \\
&\stackrel{(\ref{urops})}{=}& \langle u, \rho' \otp \Psi_*(\tau) \rangle \\
&=& \langle L_{\Psi_*(\tau)} (u), \rho' \rangle,
\end{eqnarray*}
which yields the claim.
\end{proof}
We now come to the main result of this section, which establishes a converse of the above Theorem 
for a very wide class of von Neumann algebras. It may be interpreted as a result on ``automatic normality'', i.e., 
some simple, purely algebraic condition -- namely commuting of the associated amplifications -- automatically implies the 
strong topological property of normality. Besides one 
very mild technical condition, the first von Neumann algebra, 
$\em$, can be completely arbitrary, the second one, $\en$, is only 
assumed to be properly infinite; we will show by an 
example that 
without the latter condition, the conclusion does not hold in general. 
The weak technical assumption mentioned says, roughly speaking, that the von Neumann algebras involved should not 
be acting on Hilbert spaces of pathologically large dimension. We remark that for separably acting 
von Neumann algebras $\em$ and $\en$, or, more generally, in case $\em \tt \en$ is countably decomposable, 
our assumption is trivially satisfied. 
%
In order to formulate 
the condition precisely, 
we use the following natural terminology which has been introduced in \cite{ich1}. 
\begin{defi} \label{ichhu} 
Let $\rr \subseteq \B(\H)$ be a von Neumann algebra. Then we define the $\mathrm{de\-com\-po\-sa\-bi\-lity}$ 
$\mathrm{number}$ of $\rr$, 
denoted by $\dec(\rr)$, to be the 
smallest cardinal number $\kappa$ such that 
every family of non-zero pairwise 
orthogonal projections in $\rr$ 
has at most cardinality 
$\kappa$. 
\end{defi} 
\begin{bem} 
Of course, a von Neumann algebra $\rr$ is countably decomposable if and only if $\dec(\rr) \leq \aleph_0$. 
\end{bem} 
We are now ready to state the main theorem of this section. 
\begin{satz} \label{umk} 
Let $\em$ and $\en$ be von Neumann algebras such that $\dec (\em \tt \en)$ is a non-measurable cardinal, 
and suppose $\en$ is properly infinite. If $\Phi \in \C\B(\em)$ is such that 
$$(\Phi \otp \id_\en) ~(\id_\em \otp \Psi) = (\id_\em \otp \Psi) ~(\Phi \otp \id_\en)$$
for all $\Psi \in \C\B(\en)$, then $\Phi$ must be normal. 
\end{satz} 
\begin{bem} 
The assumption that the decomposability number of $\em \tt \en$ be non-measurable, is very weak and in fact just excludes 
a set-theoretic pathology. We briefly recall that a cardinal $\kappa$ is said to be (real-valued) measurable if for every 
set $\Gamma$ of cardinality $\kappa$, there exists a diffused probability measure on the power set 
$\mathfrak{P}(\Gamma)$. -- Measurability is a property of ``large'' cardinals ($\aleph_0$ is of course not measurable). 
In order to demonstrate how weak the restriction to non-measurable cardinals is, let us just mention the fact (cf.\ \cite{kama}, 
\S 1, p.\ 106, 108) that 
the existence of measurable cardinals cannot be proved in ZFC (= the axioms of Zermelo--Fraenkel + the axiom of choice), 
and that it is consistent with ZFC to assume the non-existence of measurable cardinals (\cite{hand}, \S 4, Thm.\ 4.14, p.\ 972). 
For a further discussion, we refer to \cite{ich1}, and the references therein. 
\end{bem} 
We begin now with the preparations needed for the proof of Theorem \ref{umk}. 
We first note that, even though 
requiring a non-measurable decomposability number means only excluding exotic von Neumann algebras, it is already enough 
to ensure a very pleasant property of the predual. The latter 
is actually of 
great (technical) importance in the proof of Theorem \ref{umk} and is therefore stated explicitly in the following. 
\begin{prop} \label{mazur} 
Let $\rr$ be a von Neumann algebra. 
Then $\dec(\rr)$ is non-measurable if and only if the predual $\rr_*$ has Mazur's property (i.e., 
$w^*$-sequentially continuous functionals on $\rr$ are $w^*$-continuous -- and hence belong to $\rr_*$). 
\end{prop}
\begin{proof} 
This is 
Thm.\ 3.11 in \cite{ich1}. 
\end{proof}
The crucial idea in the proof of Theorem \ref{umk} is to use a variant of 
the concept of a \textit{countable spiral nebula}, which has been introduced in \cite{wiho}, \S 1.1 (cf.\ ibid., Lemma 1.5). 
\begin{defi} \label{nebel} 
Let $\rr$ be a von Neumann algebra. A sequence $(\al_n , e_n)_{n \in \N}$ of ${}^*$-automorphisms $\al_n$ on $\rr$ and 
projections $e_n$ in $\rr$ will be called a $\mathrm{countable ~reduced ~spiral ~nebula}$ on $\rr$ if the following holds true: 
$$e_n \leq e_{n+1} , ~{\mathrm{WOT}}-\lim_n e_n = 1 , ~ \al_m (e_m) \perp \al_n (e_n) ~\mathrm{for ~all}~ m,n \in \N_0 , 
m \not= n .$$ 
\end{defi} 
The following Proposition is the crucial step in order to establish Theorem \ref{umk}. It shows 
that a countable reduced spiral nebula may be used to derive a result on automatic normality -- and not 
only on automatic boundedness, as is done in \cite{wiho}, Lemma 1.5. 
\begin{prop} \label{hauprop} 
Let $\rr$ be a von Neumann algebra such that $\dec(\rr)$ is non-measurable, and suppose that there is 
a countable reduced spiral nebula $(\al_n , e_n)$ on $\rr$. 
\begin{itemize} 
\item[(i)] Put $l_k \dogl \frac{k (k+1)}{2}$, and let $\mathfrak{F}$ be a free ultrafilter on $\N$. Then the 
operators 
$$\Psi_n \dogl w^*-\lim_{k \rightarrow \mathfrak{F}} \al^{-1}_{l_k + n}$$
are completely positive and unital. 
\item[(ii)] If $\Phi \in \B(\rr)$ commutes with all elements of the set 
$S \dogl \overline{\{ \al_n^{-1} \mid n \in \N_0 \}}^{w^*} \subseteq \C\P(\rr)$, 
then $\Phi$ is normal. 
\end{itemize}
\end{prop} 
\begin{proof} 
The statement (i) is clear. In order to prove (ii), 
in view of Proposition \ref{mazur}, we only have to show that 
if $(x_n)_{n \in \N_0} \subseteq \Ball (\rr)$ converges $w^*$ to $0$, then the same is true for the sequence $(\Phi (x_n))_n$. 
It suffices to prove that if $(\Phi (x_{n_i}))_i$ is a $w^*$-convergent subnet, then its limit must be $0$. 
\par 
The projections $\al_n (e_n)$ being pairwise orthogonal, noting that $l_{\kappa +1} = l_\kappa + \kappa +1$, we see that 
$$\widetilde{x} \dogl \sum_{\kappa=0}^{\infty} \sum_{\nu = 0}^{\kappa} 
\al_{l_\kappa + \nu} (e_{l_\kappa + \nu} x_\nu e_{l_\kappa + \nu})  \in \rr .$$ 
Now we deduce (cf.\ the proof of \cite{wiho}, Lemma 1.5) 
that 
\begin{equation} 
\Psi_n (\widetilde{x}) = x_n ~\mathrm{for ~all}~ n \in \N_0 . \label{faktori} 
\end{equation} 
Let $\Psi \dogl w^*-\lim_{j} \Psi_{n_{i_{j}}} \in \Ball(\C\B(\rr))$ be a $w^*$-cluster point of the net $(\Psi_{n_{i}})_i$. 
Then $\Psi \in S$, and we obtain, using (\ref{faktori}): 
\begin{eqnarray*} 
\lim_i \Phi (x_{n_i}) &=& \lim_i \Phi (\Psi_{n_i} (\widetilde{x})) 
= \lim_i \Psi_{n_i} (\Phi (\widetilde{x})) \\ 
&=& \lim_j \Psi_{n_{i_{j}}} (\Phi (\widetilde{x})) 
= \Psi (\Phi (\widetilde{x})) \\ 
&=& \Phi (\Psi (\widetilde{x})) 
= \Phi (\lim_j \Psi_{n_{i_{j}}} (\widetilde{x})) \\ 
&=& \Phi (\lim_j x_{n_{i_{j}}}) 
= \Phi (0) = 0 , 
\end{eqnarray*} 
which finishes the proof. 
\end{proof} 
We are now sufficiently prepared for Theorem \ref{umk}. 
\begin{proof}
Since $\en$ is properly infinite, it is 
isomorphic to $\en \tt \B(\ell_2(\zett))$; cf., e.g., Appendix C, ``Theorem'', in \cite{vand}. By applying the argument in the first 
part of the proof of Prop.\ 2.2 in \cite{wiho}, we 
can construct a countable reduced spiral nebula $(\al_n , e_n)$ on $\B(\ell_2(\zett))$. Hence, 
$(\widetilde{\al}_n , \widetilde{e}_n) \dogl (\id_\em \otimes \id_\en \otimes \al_n , 1_\em \otimes 1_\en \otimes e_n)$ 
is a countable reduced spiral nebula on $\em \tt \en \tt \B(\ell_2(\zett)) = \em \tt \en$. 
We will now apply Proposition \ref{hauprop} to the von Neumann algebra $\em \tt \en$ and the mapping 
$\Phi \otimes \id_\en \in \C\B(\em \tt \en)$ to show that the latter, and hence $\Phi$ itself, is normal. -- To this end, 
we only have to prove that $\Phi \otimes \id_\en$ commutes with all elements of the set 
$S \dogl \overline{\{ \widetilde{\al}_n^{-1} \mid n \in \N_0 \}}^{w^*}$. But the amplification in the first variable, $\Chi_\em$, 
is $w^*$-$w^*$-continuous (this is established in the same fashion as we did for $\Chi_\em$ in Theorem \ref{amsatz}), 
and this entails that every element of $S$ is of the form $\id_\em \otimes \Psi$ for some $\Psi \in \C\B(\en)$. 
\end{proof} 
\begin{bem} 
The assumption of $\en$ being properly infinite 
may be interpreted as a certain richness condition which is necessary to be imposed on the second von Neumann algebra. 
Indeed, if in contrast $\en$ is, e.g., a finite type I factor, the statement is wrong in general. This is easily seen as follows: 
In this case, $\en = \B(\ell_2^n) = M_n$ for some $n \in \N$, and hence $\C\B(\en) = (\Te (\ell_2^n) \tp \B(\ell_2^n))^* 
= (\Te (\ell_2^n) \tp \Ka(\ell_2^n))^* = \C\B^\sigma (\en)$. 
Now let $\em$ be an arbitrary von Neumann algebra such that the singular part of the Tomiyama--Takesaki decomposition of $\C\B(\em)$ is 
non-trivial, i.e., $\C\B^{s}(\em) \not= (0)$. Then by Theorem \ref{commute}, every $\Phi \in \C\B^{s}(\em)$ will satisfy 
the commutation relation in Theorem \ref{umk} for all $\Psi \in \C\B(\en)=\C\B^\sigma (\en)$, even though it is singular. 
\end{bem}
We finish our discussion with 
considering again an arbitrary admissible pair ($\em$, $\en$). 
In view of Theorem \ref{commute}, it would be natural, following the model of the tensor product of bounded operators on
Hilbert spaces, to define a tensor product of
$\Phi \in \C\B(\em)$ and
$\Psi \in \C\B^{\sigma}(\en)$ by setting:
\begin{eqnarray*}
\Phi \overline{\otimes} \Psi \dogl (\Phi \otp \id_{\en}) ~(\id_{\em} \otp \Psi) = (\id_{\em} \otp \Psi) ~(\Phi \otp \id_{\en}),
\end{eqnarray*}
which defines an operator in $\C\B(\em \tt \en)$. It would be interesting to go even further and consider, for arbitrary
$\Phi \in \C\B(\em)$ and $\Psi \in \C\B(\en)$, the
two ``tensor products'' $\Phi \widetilde{\otimes}_1 \Psi$ and $\Phi \widetilde{\otimes}_2 \Psi$, defined respectively
through
$$\Phi \widetilde{\otimes}_1 \Psi \dogl (\Phi \otp \id_{\en}) ~(\id_{\em} \otp \Psi) \in \C\B(\em \tt \en)$$
and
$$\Phi \widetilde{\otimes}_2 \Psi \dogl (\id_{\em} \otp \Psi) ~(\Phi \otp \id_{\en}) \in \C\B(\em \tt \en).$$
By Theorem \ref{umk}, we know that these are different in general, even in the case of von Neumann algebras $\em$ and $\en$. 
We wish to close this section by stressing the resemblance of
these considerations to the well-known construction of the two Arens products on the bidual of a Banach algebra, which in general
are different. From this point of view, the subalgebra $\C\B^{\sigma}(\en)$ could play the role of what in the context
of Banach algebras is known as the topological centre (for the latter term, see, e.g., \cite{dales}, Def.\ 2.6.19). 

\section{Uniquenes of the amplification and a generalization 
of the Ge--Kadison Lemma} \label{ge} 

In the following, we 
will show that under different natural conditions, 
an algebraic amplification is uniquely determined. 
This 
leads in particular to a remarkable generalization of the so-called Ge--Kadison Lemma 
(see Proposition \ref{ichgekad}), which is obtained as an application of our Theorem \ref{commute}. 
\par 
But before establishing positive results, we briefly point out why 
in general, an algebraic amplification 
is highly non-unique. In \cite{tom}, p.\ 28, Tomiyama states, in the context of von Neumann 
algebras, that the ``product projection'' of two 
projections of norm one ``might not be unique'' (of course, in our case, the second projection of norm one 
would be the identity mapping). The following simple argument shows that whenever $\em$ and $\en$ 
are infinite-dimensional dual operator spaces, there are uncountably many different algebraic amplifications 
of a map $\Phi \in \C\B(\em)$ to a completely bounded map on $\em \tt \en$, regardless of $\Phi$ being normal or not. 
Namely, for any non-zero functional 
$\varphi \in ( \em \tt \en )^*$ which vanishes on $\em \ti \en$, and any non-zero vector $v \in \em \tt \en$, 
$$\Chi_{\en}^{\varphi, v} (\Phi) \dogl \Chi_{\en} (\Phi) - \langle \varphi , \Chi_{\en} (\Phi) ( \cdot ) \rangle ~v$$ 
defines such an algebraic amplification. 
Furthermore, the amplification $\Chi_{\en}^{\varphi, v}$ can of course be taken arbitrarily close to an isometry: 
for every $\varepsilon > 0$, an obvious choice of $\varphi$ and $v$ yields an algebraic amplification 
$\Chi_{\en}^{\varphi, v} : \C\B (\em) \longrightarrow \C\B ( \em \tt \en )$ such that 
$$(1-\varepsilon ) ~\| \Phi \|_{\mathrm{cb}} \leq \| \Chi_{\en}^{\varphi, v} (\Phi) \|_{\mathrm{cb}} 
\leq (1+\varepsilon ) ~\| \Phi \|_{\mathrm{cb}}$$ 
for all $\Phi \in \C\B(\em)$. 
\par 
We 
begin our investigation by noting 
a criterion which establishes the uniqueness of amplification 
for a large class of
von Neumann algebras, 
though only assuming part of the requirements listed in Definition \ref{amdef}. 
\begin{prop} \label{eindam}
Let $\em$ and $\en$ be von Neumann algebras, where $\em$ is an injective factor.
If $\Chi': \C\B(\em) \longrightarrow \C\B(\em \tt \en)$ is a bounded linear mapping 
which satisfies the algebraic amplification condition (AAC) and properties
(iii) and (iv) in Definition \ref{amdef}, then $\Chi' = \Chi_\en$.
\end{prop}
\begin{proof}
Owing to \cite{chsm}, Thm.\ 4.2, we have the density relation:
$${\overline{\C\B^{\sigma}(\em)}}^{w^*} = \C\B(\em).$$
Hence, due to 
their $w^*$-$w^*$-continuity (condition (iii)), $\Chi_\en$ and $\Chi'$ 
will be identical if they only coincide
on $\C\B^{\sigma}(\em)$. But if $\Phi \in \C\B^{\sigma}(\em)$, by condition (iv), the mappings
$\Chi_\en (\Phi)$ and $\Chi'(\Phi)$
are normal, so that
to establish their equality, it suffices to show that they coincide on all elementary tensors $S \otp T \in \em \tt \en$.
But this is true since they both fulfill condition 
(AAC). 
\end{proof} 
%
We will now present our generalization of the Ge--Kadison Lemma which 
was first obtained 
in \cite{gekad} 
(ibid., Lemma F) and used in the 
solution of the splitting problem for tensor products of von Neumann algebras (in the factor case). 
Furtheron, a version of the Lemma was again used by Str\u{a}til\u{a} and Zsid\'{o} in \cite{stzs} in order to 
derive an extremely general commutation theorem which extends both Tomita's classical 
commutant theorem and the above mentioned splitting result (\cite{stzs}, Thm.\ 4.7). 
The result we give below generalizes the Lemma in the version as stated (and proved) in \cite{stzs}, \S 3.4, which 
the authors refer to as a ``smart technical device''. In \cite{stzs}, it is applied to \textit{normal} 
conditional expectations of von Neumann algebras (ibid., \S 3.5); of course, these are in particular 
\textit{completely positive}. Our generalization yields an analogous result 
which in turn can be applied 
to \textit{arbitrary} completely bounded mappings. This fact may be of interest in the further development 
of the subject treated in \cite{gekad} and \cite{stzs}. Furthermore, 
we prove that the Lemma not only holds for 
von Neumann algebras 
but for arbitrary admissible pairs, a situation which is not considered in 
the latter articles. Hence, our version is likely to prove useful in the interesting problem 
of obtaining a splitting theorem 
in the more general context of $w^*$-spatial tensor products of 
ultraweakly closed subspaces -- a question which is actually hinted at in 
\cite{gekad} (p.\ 457). 
%
\begin{prop} \label{ichgekad} 
Let $(\em , \en)$ be an admissible pair, and $\Phi \in \C\B(\em)$ an arbitrary completely bounded operator. 
Suppose $\Theta : \em \tt \en \longrightarrow \em \tt \en$ is $\mathrm{any}$ map 
which satisfies, for some non-zero element $n \in \en$: 
\begin{itemize} 
\item[(i)] $\Theta$ commutes with the 
slice maps $\id_\em \otimes \tau n$ ($\tau \in \en_*$) 
\item[(ii)] $\Theta$ coincides with $\Phi \otimes \id_\en$ on $\em \otimes n$. 
\end{itemize} 
Then we must have $\Theta = \Phi \otimes \id_\en$. 
\par 
In particular, if $\Phi \in \C\B(\em)$ and $(\Phi_\al)_\al \subseteq \C\B(\em)$ is a bounded net converging $w^*$ to $\Phi$, 
then 
$\Phi_\al \otimes \id_\en \stackrel{w^*}{\longrightarrow} \Phi \otimes \id_\en$. 
\end{prop}
\begin{bem} 
(i) The known version of the Lemma supposes that $\em$ and $\en$ are von Neumann algebras, and that 
the mapping $\Phi$ is normal and completely positive. As shown by the above, all these 
assumptions can be either weakened or even just dropped. 
\par 
(ii) The last statement in the Proposition is even true $\mathrm{without}$ assuming the net $(\Phi_\al)_\al$ 
to be bounded, and in this generality follows directly from 
our Theorem \ref{amsatz} (property (iii) in Definition \ref{amdef}). Nevertheless, this ``unbounded'' version 
does not appear in the literature. (We remark that the assumption of boundedness is missing in the statement of Lemma 3.4 in 
\cite{stzs}, though needed in the proof presented there.) 
We shall 
point out below 
how the ``bounded'' version can be easily obtained by using the first part of the Proposition (cf.\ \cite{stzs}, \S 3.4). 
\end{bem} 
\begin{proof} 
Using the first result of the last section, our argument follows the same lines as the one given in \cite{stzs}; we 
present the proof because of its brevity. -- Fix a non-zero element $m \in \em$. We obtain for every 
$u \in \em \tt \en$, $\rho \in \em_*$ and $\tau \in \en_*$: 
\begin{eqnarray*} 
\langle \Theta (u), \rho \otimes \tau \rangle ~m \otimes n &=& 
\left[ (\rho m \otimes \id_\en ) ~( \id_\em \otimes \tau n ) ~\Theta  \right] (u) \\ 
&=& \left[ (\rho m \otimes \id_\en ) ~\Theta ~( \id_\em \otimes \tau n )  \right] (u) \quad \mathrm{by ~(i)} \\ 
&=& \left[ (\rho m \otimes \id_\en ) ~(\Phi \otimes \id_\en ) ~( \id_\em \otimes \tau n )  \right] (u) \quad \mathrm{by ~(ii)} \\ 
&=& \left[ (\rho m \otimes \id_\en ) ~( \id_\em \otimes \tau n ) ~(\Phi \otimes \id_\en ) \right] (u) 
\quad \mathrm{by ~Theorem ~\ref{commute}} \\ 
&=& \langle (\Phi \otimes \id_\en) (u), \rho \otimes \tau \rangle ~m \otimes n , 
\end{eqnarray*} 
which proves the first statement in the Proposition. 
In order to prove the second assertion, we only have to show that if $\Theta$ is any $w^*$-cluster point of the net 
$(\Phi_\al \otimes \id_\en)_\al$, then $\Theta = \Phi \otimes \id_\en$. But $\Theta$ trivially satisfies condition (ii) 
in our Proposition, and using the normality of the mappings $\id_\em \otimes \tau n$, 
another application of Theorem \ref{commute} shows that $\Theta$ also meets condition (i), whence the first part of the 
Proposition finishes the proof. 
\end{proof}

\section{Completely bounded module homomorphisms and Tomiyama's slice maps} \label{sechs}
The main objective of \cite{mnw} was to establish the following remarkable property of the amplification
mapping
$$\C\B(\B(\H)) \ni \Phi \mapsto \Phi^{(\infty)} \in \C\B(\B(\H) \tt \B(\K)).$$
If $\rr \subseteq \B(\H)$ is a von Neumann subalgebra, and if $\Phi \in \C\B_{\rr}(\B(\H))$ is an
$\rr$-bimodule homomorphism, then $\Phi^{(\infty)}$ is an $\rr \tt \B(\K)$-bimodule homomorphism, i.e.,
$\Phi^{(\infty)} \in \C\B_{\rr \tt \B(\K)}(\B(\H) \tt \B(\K))$; see \cite{mnw}, Prop.\ 2.2 (and
\cite{bim}, Thm.\ 4.2, for a
generalization).
\par
An important special
situation, which nevertheless shows the 
strength and encodes the essential statement of the result, is the case where 
$\rr = {\mathbb{C}} ~\1_{\B(\H)}$. Then the above theorem reads as follows:
\par
$(**) ~{}$ Every $\Phi \in \C\B(\B(\H))$
can be amplified to an
$\1_{\B(\H)} \otimes \B(\K)$-bimodule homomorphism $\Phi^{(\infty)} \in \C\B(\B(\H) \tt \B(\K))$.
\par
The aim of this section is to prove
an analogous
result where now $\B(\H)$ and $\B(\K)$ are replaced by arbitrary
von Neumann algebras $\em$ and $\en$, using our amplification $\Chi_\en$. In this section, $\em$ and $\en$ will
always denote von Neumann algebras.
\par
The proofs of \cite{mnw}, Prop.\ 2.2 and \cite{bim}, Thm.\ 4.2
use either the Wittstock decomposition theorem or
a corresponding result by Paulsen (\cite{paul}, Lemma 7.1) to reduce to the case of a completely positive mapping,
where the assertion is then obtained by a fairly subtle analysis involving the Schwarz inequality.
In contrast to this,
the explicit form of $\Chi_\en$ provides us with a
direct route to
the above announced goal; indeed, as we shall see,
the
statement asserting the module homomorphism property of the amplified mapping
reduces to a corresponding statement about Tomiyama's slice maps which in turn can be verified in an elementary fashion (see
Lemma \ref{easy} below).
\par
In order to
prove our result, a little preparation is needed.
In the
sequel, we restrict ourselves to left slice maps; analogous statements hold of course in the case of
right slice maps. -- As is well-known (and easy to see), Tomiyama's left slice maps satisfy the following module
homomorphism property:
$$L_\tau ( (a \otimes \1_\en) u (b \otimes \1_\en) ) = a L_\tau (u) b,$$
where $u \in \em \tt \en$, $a, b \in \em$, $\tau \in \en_*$.
In the
following (elementary) lemma, we shall
establish an analogous property involving elementary tensors of the form $\1_\em \otimes a$ and $\1_\em \otimes b$ for
$a, b \in \en$. To this end, we briefly recall that for every von Neumann algebra $\en$, the predual $\en_*$
is an $\en$-bimodule in a very natural manner, the actions being given by
$$\langle \tau \cdot a, b \rangle = \langle \tau, ab \rangle \quad \mathrm{and} \quad
\langle a \cdot \tau , b \rangle = \langle \tau, ba \rangle,$$
for $a, b \in \en$, $\tau \in \en_*$.
For later purposes, we note the (easily verified) equation:
\begin{eqnarray}
(\1_\em \otimes a) \cdot (\rho \otimes \tau) \cdot (\1_\em \otimes b) = \rho \otimes (a \cdot \tau \cdot b), \label{gleichun}
\end{eqnarray}
where $a, b \in \en$, $\rho \in \em_*$, $\tau \in \en_*$.
We now come to
\begin{lemma} \label{easy}
Let $\em$ and $\en$ be von Neumann algebras. Then Tomiyama's left slice map satisfies:
$$L_\tau ( (\1_\em \otimes a) u (\1_\em \otimes b) ) = L_{b \cdot \tau \cdot a} (u),$$
where $u \in \em \tt \en$, $a, b \in \en$, $\tau \in \en_*$.
\end{lemma}
\begin{proof}
Due to the normality of
left slice maps, it suffices to prove the statement
for elementary tensors $u=S \otimes T \in \em \tt \en$.
But in this case, we obtain:
\begin{eqnarray*}
L_\tau ( (\1_\em \otimes a) u (\1_\em \otimes b) ) &=& L_\tau (S \otimes aTb) \\
&=& \langle \tau, aTb \rangle ~S \\
&=& \langle b \cdot \tau \cdot a, T \rangle ~S \\
&=& L_{b \cdot \tau \cdot a} (u),
\end{eqnarray*}
whence the desired equality follows.
\end{proof}
The above
statement
will now be used to derive a corresponding module homomorphism property shared by the
amplification of a completely bounded mapping, as announced at the beginning of this section.
%
\begin{satz} \label{modulhom}
Let $\em$ and $\en$ be von Neumann algebras. Then for every $\Phi \in \C\B(\em)$, the
amplified mapping $\Chi_\en (\Phi) \in \C\B(\em \tt \en)$ is an
$\1_{\em} \otimes \en$-bimodule homomorphism.
\end{satz}
\begin{proof}
As in the proof of Theorem \ref{amsatz}, we denote by $\kappa$ the pre-adjoint mapping of $\Chi_\en$. Fix
$u \in \em \tt \en$, $a, b \in \en$, $\rho \in \em_*$ and $\tau \in \en_*$.
We write $a' \dogl \1_{\em} \otimes a$, $b' \dogl \1_{\em} \otimes b$.
\par
We first note that
\begin{eqnarray}
\kappa (\rho \otimes \tau \otimes (a' u b')) =
\kappa (\rho \otimes (b \cdot \tau \cdot a) \otimes u). \label{gleichungg}
\end{eqnarray}
For
by equation (\ref{chi}), in order to establish (\ref{gleichungg}),
we only have to show that $$L_\tau (a' u b') = L_{b \cdot \tau \cdot a}(u),$$ which
in turn is precisely the statement of Lemma \ref{easy}.
\par
Now we
obtain for $\Phi \in \C\B(\em)$:
\begin{eqnarray*}
\langle \Chi_\en(\Phi) (a' u b'), \rho \otimes \tau \rangle &=& \langle \Phi, \kappa(\rho \otimes \tau \otimes (a' u b')) \rangle \\
&\stackrel{(\ref{gleichungg})}{=}& \langle \Phi, \kappa (\rho \otimes (b \cdot \tau \cdot a) \otimes u) \rangle \\
&\stackrel{(\ref{gleichun})}{=}& \langle \Phi, \kappa [ (b' \cdot (\rho \otimes \tau) \cdot a') \otimes u] \rangle \\
&=& \langle \Chi_\en(\Phi) (u), b' \cdot (\rho \otimes \tau) \cdot a' \rangle \\
&=& \langle a' \Chi_\en(\Phi) (u) b', \rho \otimes \tau \rangle,
\end{eqnarray*}
which
finishes the proof.
\end{proof}
We
finish by
remarking that Theorem \ref{modulhom} yields in particular the statement $(**)$ given above; for if
$\Phi \in \C\B(\B(\H))$, we have $\Chi_{\B(\K)}(\Phi)=\Phi^{(\infty)}$, as noted in section \ref{disc}, equation (\ref{glinf}). 
%
%
%

\vspace{0.7cm} 
\noindent
{\large{\textbf{Acknowledgments}}}
\\
\noindent
The author is indebted to Prof.\ Dr.\ G.\ Wittstock for many
very fruitful discussions on the subject of the present paper.

\vspace{0.2cm}
\noindent
{\sc{Author's address:}}
\\
{\textit{
Department of Mathematical Sciences\\
University of Alberta\\
Edmonton, Alberta\\
Canada T6G 2G1\\
E-mail:}}
mneufang@math.ualberta.ca
\end{document}